\documentclass[preprint]{elsarticle}

\usepackage{graphicx}
\usepackage{geometry}\usepackage{amssymb,bm}
\usepackage{amsmath}
\usepackage{subcaption}
\usepackage{hyperref}
\usepackage{xcolor}
\usepackage{soul}
\usepackage{appendix}
\usepackage{algorithm}
\usepackage[noend]{algpseudocode}

\usepackage{lineno}

\DeclareMathAlphabet\mathbfcal{OMS}{cmsy}{b}{n}
\newcommand{\RM}{\bm{\Lambda}}
\newcommand{\RMI}{\bm{\Lambda}_0}
\newcommand{\RMT}{\bm{\Lambda}_t}
\newcommand{\RV}{\bm{\psi}}
\newcommand{\magRV}{\psi}
\newcommand{\diag}{\rm diag}

\title{A cell-centred finite volume formulation of geometrically-exact Simo-Reissner beams with arbitrary initial curvatures}

\author[1,5,6,7,9]{Seevani Bali}

\author[2,6]{\v{Z}eljko Tukovi\'{c}}

\author[1,4,6,9]{Philip Cardiff\corref{cor1}\fnref{fn1}}

\author[3,4,6,8,9]{Alojz Ivankovi\'{c}}

\author[1,5,6,7,9]{Vikram Pakrashi}

\cortext[cor1]{Corresponding author}
\fntext[fn1]{philip.cardiff@ucd.ie}

\address[1]{School of Mechanical and Materials Engineering, University College Dublin, Ireland}
\address[2]{Faculty of Mechanical Engineering and Naval Architecture, University of Zagreb, Croatia}
\address[3]{Professor, School of Mechanical and Materials Engineering, University College Dublin, Ireland}
\address[4]{SFI I-Form Centre, University College Dublin, Ireland}
\address[5]{SFI MaREI Centre, University College Dublin, Ireland}
\address[6]{Bekaert University Technology Centre, School of Mechanical and Materials Engineering, University College Dublin, Ireland}
\address[7]{Dynamical Systems and Risk Laboratory, School of Mechanical and Materials Engineering, University College Dublin, Ireland}
\address[8]{UCD Centre of Adhesion and Adhesives, University College Dublin, Ireland}
\address[9]{UCD Centre for Mechanics, University College Dublin, Ireland}
\begin{document}

\begin{abstract}
This paper presents a novel total Lagrangian cell-centred finite volume formulation of geometrically exact beams with arbitrary initial curvature undergoing large displacements and finite rotations. The choice of rotation parametrisation, the mathematical formulation of the beam kinematics, conjugate strain measures and the linearisation of the strong form of governing equations is described. The finite volume based discretisation of the computational domain and the governing equations for each computational volume are presented. The discretised integral form of the equilibrium equations are solved using a block-coupled Newton-Raphson solution procedure. The efficacy of the proposed methodology is presented by comparing the simulated numerical results with classic benchmark test cases available in the literature. The objectivity of strain measures for the current formulation and mesh convergence studies for both initially straight and curved beam configurations are also discussed.


\end{abstract}



\maketitle

\section{INTRODUCTION}
\label{sec:Intro}

\noindent The mathematical modelling of nonlinear beams has flourished in the past few decades with their applicability spanning various fields of engineering. Meier et al. ~\cite{meier2019geometrically} presented a concise review of the existing beam theories. These formulations, in particular, and the numerical computations of solid mechanics problems, in general, are mostly analysed using the finite element (FE) approach; however, since the $1980$s, the applicability of the simple and conservative finite volume (FV) methods to solid mechanics problems has become increasingly popular and is evolving rapidly~\cite{cardiff2021thirty}. Particularly in the purview of beams, the FV implementation of Euler-Bernoulli and shear-deformable Timoshenko beam theories for both straight and curved beam configurations, can be found in \cite{moosavi2011orthogonal,fallah2006displacement,isic2007comparison,fallah2014finite,fallah2018displacement,jing2016static}.
Additionally, these works have been extended to study the stability and buckling analysis of beam-columns as well as plates \cite{fallah2006extension,fallah2013finite,isic2007numerical,fallah2004cell,fallah2005new,fallah2006finite,wheel1997finite}. In contrast to the FE implementation,  Fallah and Ghanbari \cite{fallah2018displacement} and Wheel \cite{wheel1997finite} reported that the FV formulation did not suffer from ``shear-locking" for both thin Timoshenko beams and thin Mindlin-Reissner plate analysis respectively. The geometrically exact Simo-Reissner beam formulation is the most general nonlinear 3-D beam theory capable of dealing with finite displacements and rotations. The FV formulation of geometrically exact beams was first investigated by Tukovi\'{c} et al. \cite{tukovic2020afinite}. The current work builds on the former developed formulations for quasi-static shear-deformable geometrically exact beams and provides a lucid comparison between the FE based classic benchmark cases and the proposed FV methodology. This is the first article to present a FV formulation for beams subjected to finite displacements and rotations. 


The formulations of a beam in a geometrically exact sense considering 3-D rotations and its FE approximation, was first proposed by Simo et al.~\cite{Simo198555,simo1986three}. The involvement of finite rotations, elements of the nonlinear differentiable manifold $SO(3)$ and the complexity of interpolating large rotations~\cite{argyris1982excursion,Ibrahimbegovic199749} led to three major techniques of addressing them, viz., (a) incremental and total rotational vector based parametrisations~\cite{makinen2007total,ibrahimbegovic1995computational,jelenic1995kinematically,Marino2016383}, (b) quaternion based rotation interpolation~\cite{ghosh2008consistent,zupan2009quaternion,geradin1988kinematics,kapania2003formulation} and, (c) co-rotational beam formulations~\cite{crisfield1990consistent,cardona1988beam}, to name a few. Romero \cite{romero2004interpolation} presents a concise review of different rotation interpolations used for geometrically exact beams. In the current work, the finite rotations are parameterised using rotational vectors and the rotational strain measures are updated using incremental rotation vectors from the previously converged configuration of the beam. The loss of objectivity in strain measures for FE approximation of interpolated rotations was first pointed out be Jeleni\'{c} et al. \cite{Crisfield19991125} which was followed by strain invariant formulations of the geometrically exact beams \cite{jelenic1999geometrically,romero2002objective,zupan2003finite,ghosh2009frame}. For the current work, the objectivity of adopted strain measures are numerically verified using a test case for the current work (Section \ref{subsec:objectivity}). As a recent contribution, Meier et al. \cite{meier2014objective} developed geometrically exact Kirchhoff-Love formulation for slender rod geometries. For a general review on other research advances in this field, like discrete Cosserat rod kinematics, implementation of different time-stepping schemes in beam multi-body dynamics, inclusion of nonlinear constitutive laws and enhanced kinematics of beams, the interested readers are referred to articles by Meier et al. \cite{meier2019geometrically} and Chadha et al. \cite{chadha2020mathematical}.

The article is structured as follows: Section \ref{sec:Math} outlines the mathematical model in total Lagrangian form. The FV discretisation of the mathematical model is described in Section \ref{sec:numerics}. In Section \ref{sec:testCases}, the proposed methodology is evaluated on five complementary benchmark test cases, where predictions are compared to analytical solutions and existing FE benchmarks.


\section{MATHEMATICAL MODEL}
\label{sec:Math}
In this section, the mathematical formulation for a 3-D quasi-static, shear-deformable geometrically-exact Simo-Reissner beam is summarised. Subsequently, the governing equilibrium equations of spatial forces and moments, the constitutive relations and their linearised formulation are described.


\subsection{Kinematic description}\label{subsec:kinematics}

\begin{figure}[h]
	\centering
	\includegraphics[scale = 0.5]{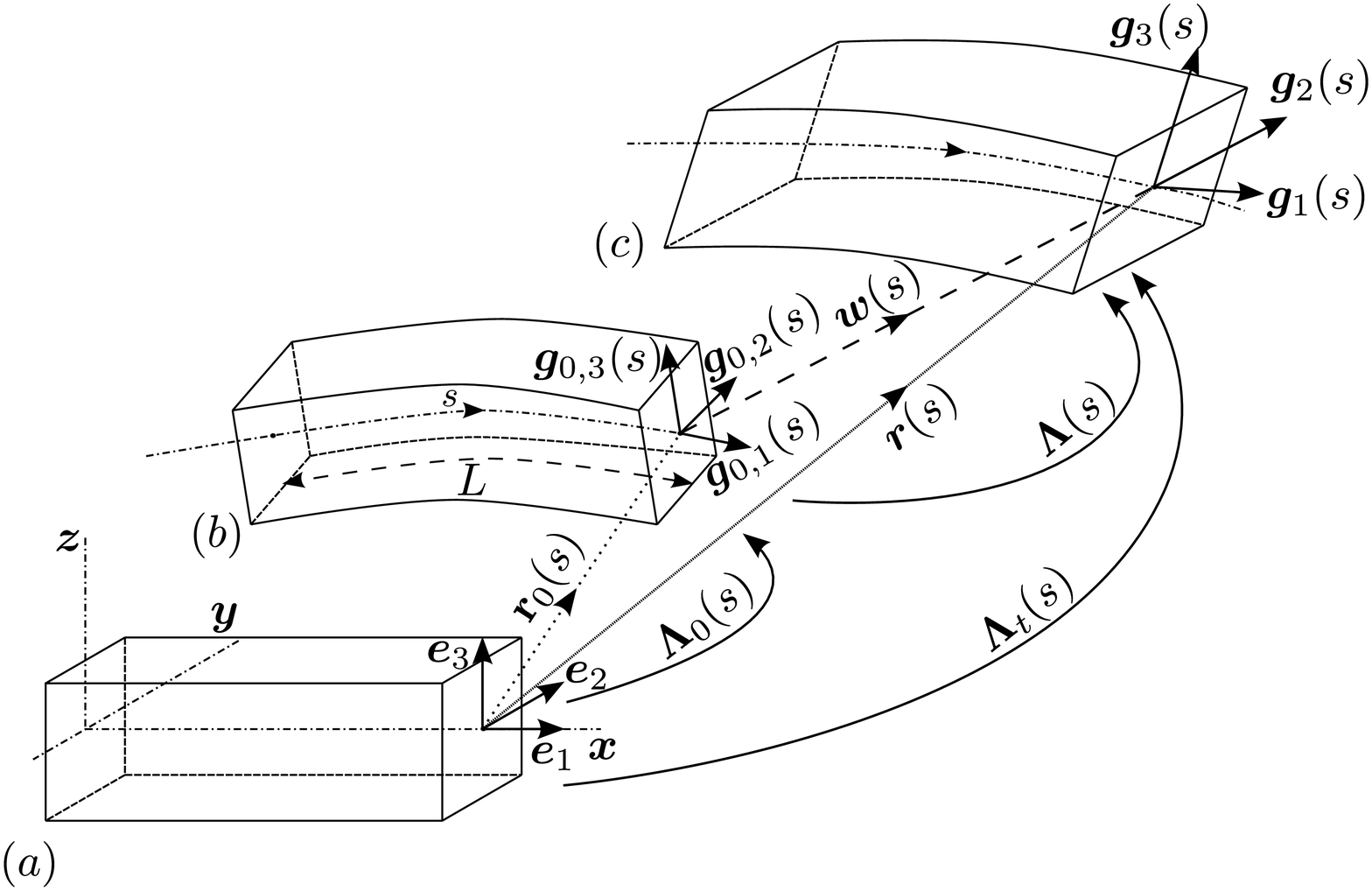}
	\caption{Beam Kinematics: (a) reference, (b) initial and, (c) deformed configurations of the beam respectively}
	\label{fig:Beam}
\end{figure}

\noindent  A total Lagrangian formulation is adopted to describe the large deformations of the beam model. Accordingly, a right-handed fixed, reference Cartesian frame defined by the orthonormal basis vectors $\bm{e}_1 = [1, 0, 0]^{\rm T}$, $\bm{e}_2 = [0, 1, 0]^{\rm T}$ and $\bm{e}_3 = [0, 0, 1]^{\rm T}$ is specified (Fig. \ref{fig:Beam}). In the reference (i.e. material) configuration, the mean line of a beam (i.e. line of centroids of the beam cross-sections) is straight and parallel to the basis vector $\bm{e}_1$ and the cross-sections of the beam are orthogonal to the basis vector $\bm{e}_1$. The principal axis of inertia of the beam cross-sections are directed along the basis vectors $\bm{e}_2$ and $\bm{e}_3$.

An initial (stress-free) configuration of the beam mean line is defined by the space curve $\bm{r}_0(s)$, where  $s \in [0,L]$ is the beam-length and $L$ is the initial length of the beam (Fig. \ref{fig:Beam}). To describe the continuous kinematic description of the moving mean line curve of the beam, a body-attached right handed orthonormal base vectors $\bm{g}_{0,1}(s)$, $\bm{g}_{0,2}(s)$ and $\bm{g}_{0,3}(s)$ are defined by the Frenet-Serret formulae, where the base vector $\bm{g}_{0,1}(s)$ is directed along the initial mean line so that

\begin{equation}
	\bm{g}_{0,1}(s) = \bm{r}_0^\prime(s)
	\label{eq:initialTangent}
\end{equation}

\noindent and the base vectors $\bm{g}_{0,2}(s)$ and $\bm{g}_{0,3}(s)$ are directed along the principal axis of inertia of the cross-section at $s$. The prime $(\cdot)^\prime$ operator in Eq. \ref{eq:initialTangent} and hence forward wherever used, denotes a derivative with respect to arc-length parameter $s$, i.e. $(\cdot)^\prime \equiv \frac{\partial (\cdot)}{\partial s}$. For an initially straight beam, the base vectors $\bm{g}_{0,i}$ and the reference bases $\bm{e}_i$ coincide, and the $\bm{g}_{0,1}(s) \equiv \bm{e}_1(s)$.

\noindent The orthonormal basis $\bm{g}_{0,1}(s)$, $\bm{g}_{0,2}(s)$, $\bm{g}_{0,3}(s)$ of the local frame and the orthonormal bases $\bm{e}_1$, $\bm{e}_2$, $\bm{e}_3$ of the global Cartesian frame are related through a linear transformation $\bm{\Lambda}_0(s)$ as
\begin{equation}
	\bm{g}_{0,i}(s) = \bm{\Lambda}_0(s)\bm{e}_{i},\quad i = 1,2,3 
\end{equation}
where $\bm{\Lambda}_0(s) \in SO(3)$ is the initial two-point second-order orthogonal rotational tensor field that defines the orientations of beam cross-sections w.r.t the reference basis. Hence, the initial configuration of the beam is fully defined by the position vector $\bm{r}_0(s)$ of the beam mean line and orientation of the cross-section at $s$ via the orthogonal rotation tensor $\bm{\Lambda}_0(s)$.

The deformed configuration of the beam mean line is defined by a space curve $\bm{r}(s)$ and the orientation of the cross-sections via the moving spatial bases $\bm{g}_{1}(s)$, $\bm{g}_{2}(s)$, $\bm{g}_{3}(s)$. Contrary to the initial configuration, the basis vector $\bm{g_{1}(s)}$ need not be directed along the deformed mean line since the beam model is capable of representing shear deformations. The base vectors $\bm{g}_{2}(s)$ and $\bm{g}_{3}(s)$, however, are still directed along the principal axes of the cross-section at $s$ and the orthonormal bases $\bm{g}_{1}(s)$, $\bm{g}_{2}(s)$, $\bm{g}_{3}(s)$ are related to $\bm{g}_{0,1}(s)$, $\bm{g}_{0,2}(s)$, $\bm{g}_{0,3}(s)$ bases by the linear transformation
\begin{equation}
	\bm{g}_{i}(s) = \bm{\Lambda}(s)\bm{g}_{0,i}(s) = 
	\bm{\Lambda}(s)\bm{\Lambda}_0(s)\bm{e}_{i}, \quad i=1,2,3 
\end{equation}
\noindent where $\bm{\Lambda}(s) \in SO(3)$ is the relative rotational matrix which rotates beam cross-sections from the initial to the deformed configuration. Additionally, the deformed mean line of the beam can also be defined using the mean line displacement vector $\bm{w}(s)$ and the mean position vector of the initial configuration by, 
\begin{equation}\label{eq:meanDefLine}
	\bm{r}(s) = \bm{r}_0(s) + \bm{w}(s)
\end{equation}
\noindent The position vector of the deformed mean line $\bm{r}(s)$ and the orientation of the orthonormal frame $\bm{\Lambda}(s)$ attached to the cross-section at $s$, ($\bm{r}(s), \bm{\Lambda}(s)$) fully define the deformed configuration of the beam. Orientation of the beam cross-sections in the deformed configuration can also be defined w.r.t to the fixed reference bases by the total rotation matrix,
\begin{equation}
	\RMT(s) = \RM(s)\RMI(s)
\end{equation}


\subsection{Balance equations and strain measures}\label{subsec:baleqAndStrain}
\noindent Assuming that the deformation from the initial to the deformed configuration is caused by distributed external forces and torques $\bm{f}$ and $\bm{t}$ respectively, the strong differential form of the balance equations are given as \cite{Simo198555},
\begin{equation}
	\bm{n}^\prime \ + \  \bm{f} = 0	
	\label{eq:forceBalance}
\end{equation}
\begin{equation}
	\bm{m}^\prime \ + \ \bm{r}^\prime\times\bm{n}  \ + \ \bm{t} = 0
	\label{eq:momentBalance}
\end{equation}
where $\bm{n}$ and $\bm{m}$ are the vectors of spatial internal forces and moments acting over the cross-section at $s$ and ``$\times$" symbol denotes the cross product between two vectors. The corresponding material counterparts ($\bm{N}$ and $\bm{M}$) are related to $\bm{n}$ and $\bm{m}$ via the pull-back mapping $\RMT^{\rm T}$ as $\bm{N} = \RMT^{\rm T}\bm{n}$ and $\bm{M} = \RMT^{\rm T}\bm{m}$~\cite{Simo198555}.

\noindent The equivalent strong integral form of the above balance equations over a length, $L$ can be expressed as,

\begin{equation}\label{eq:intForce}
	\int_{L} \bm{n}^\prime \textrm{d}L \ + \ \int_{L} \bm{f} \textrm{d}L = 0
\end{equation}

\begin{equation}\label{eq:intMoment}
	\int_{L} \bm{m}^\prime \textrm{d}L \ + \ \int_{L} (\bm{r}^\prime \times \bm{n}) \textrm{d} L \ + \ \int_{L} \bm{t} \textrm{d}L = 0
\end{equation}
%

\noindent The required strain measures follow from a geometrically exact beam theory, where the relationships between the beam configuration and the strain measures are consistent with the virtual work principle and the equilibrium equations at a deformed state regardless of the magnitude of displacements, rotations and strains \cite{simo1986three}. To that end, the equivalent strain-configuration relationships involve three translational strains and a skew-symmetric tensor of the three rotational strains \cite{Simo198555,ibrahimbegovic1995computational} (from here on, the argument $(s)$ is dropped from the terms for clarity):
\begin{equation} \label{eq:gamma}
	\bm{\Gamma} = \RMT^{\rm T}\bm{r}^\prime - \RMT^{\rm T} \bm{g}_1
\end{equation}
\begin{equation}
	\hat{\bm{K}} = \RMT^{\rm T}\RMT^\prime - \RMI^{\rm T}\RMI^\prime 
	= \RMT^{\rm T}\RMT^\prime
	\label{eq:rotStrains}
\end{equation}
where $\bm{\Gamma}$ and $\bm{K}$ are (material) translational and rotational strain measures respectively. Physically, the strain vector $\bm{\Gamma}$ represents axial tension (first entry) and shear deformation whereas the (material) rotational strain vector $\bm{K}$ represents torsion (first entry) and bending deformation of the beam body. These strain measures are the energy-conjugate pairs to the stress resultants ($\bm{N}$ and $\bm{M}$).

\noindent The hat $\hat{(\cdot)}$ operator denotes a skew-symmetric matrix associated with the corresponding (pseudo)-vector 
given by the relation, $\bm{a} \times \bm{h} = \hat{\bm{a}} \bm{h}$, $\forall$ $\bm{h} \in \mathbb{R}^3 $. In this work, the relative rotation matrix $\RM$ is parameterised in terms of its rotational vector $\RV$ as
\begin{equation} \label{eq:rodrigues}
	\RM(\RV) = \exp (\hat{\bm{\psi}}) = \mathbf{I} \
	+ \ \frac{\sin \magRV}{\magRV}\hat{\RV} \
	+ \ \frac{1-\cos\magRV}{\magRV^2} \hat{\RV}\hat{\RV}
\end{equation}
where $\mathbf{I}$ is a $3 \times 3$ identity matrix and $\magRV$ is the magnitude of rotation vector $\RV$. The $\hat{\bm{\psi}} \in \mathfrak{so(3)}$ is the skew-symmetric tensor living in the tangent space of $SO(3)$ at $\RM$ and its exponentiation yields the finite rotation $\RM \in SO(3)$. For a rotation vector based parametrisation, an alternative expression for evaluating incremental (material) rotational strain vector $\Delta \bm{K}$, as demonstrated in \cite{Ibrahimbegovic199749,ibrahimbegovic1995computational,jelenic1999geometrically} is given by,
\begin{equation} \label{eq:K}
	\Delta \bm{K} = \RMT^{\rm T} \bm{T}^{\rm T}\RV^\prime
\end{equation}
where the tangent operator $\bm{T} \in \mathfrak{so(3)}$ is defined as (``$\otimes$" denotes the dyadic product between two vectors),
\begin{equation} \label{eq:tangent}
	\bm{T} = 
		\frac{\sin \magRV}{\magRV}\mathbf{I} \
	   +\ \frac{1}{\magRV^2}\left(1-\frac{\sin\magRV}{\magRV}\right)\RV \otimes \RV \
	   +\ \frac{1-\cos\magRV}{\magRV^2}\hat{\RV} \
\end{equation}

\subsection{Constitutive relations}

\noindent The present study is limited to linear hyperelastic materials whose length-specific stored energy function is given by,

\begin{equation}\label{eq:intEnergy}
	\tilde{\Pi}_{\text{int}} = \frac{1}{2} \bm{\Gamma}^{\rm T} \bm{C}_{\rm N} \bm{\Gamma} + \frac{1}{2} \bm{K}^{\rm T} \bm{C}_{\rm M} \bm{K} 
\end{equation}

\noindent For this hyperelastic model, the material internal forces $\bm{N}$ and moments $\bm{M}$ are linearly related to the material strain measures ($\bm{\Gamma}$ and $\bm{K}$) as, $\bm{N} = \bm{C}_{\rm N}\bm{\Gamma}$ and $\bm{M} = \bm{C}_{\rm M}\bm{K}$ where $\bm{C}_{\rm N} = \diag{[EA, GA_2, GA_3]}$ and $\bm{C}_{\rm M} = \diag{[GJ, EI_2, EI_3]}$ are constant diagonal constitutive matrices. Here, $E$ and $G$ denote the elastic and shear moduli of the material respectively, and $A_{\rm p}$ and $I_{\rm p}$ ($\rm p = 1,2,3$), are the effective (material) areas and area moments about the principal axes of inertia respectively ($J \equiv I_1$). 

\subsection{Linearisation of balance equations} \label{subsec:linearisation}

\noindent To apply the Newton-Raphson solution procedure, the spatial internal forces and moments acting on the beam cross-section must be linearised. In this paper, the procedure described by Marino \cite{Marino2016383} and the spatial internal forces and moments,
%
\begin{equation}\label{eq:n}
	\bm{n} = \RMT\bm{C}_{\rm N}\bm{\Gamma}
\end{equation}
\begin{equation}\label{eq:m}
	\bm{m} = \RMT\bm{C}_{\rm M}\bm{K}
\end{equation}
are linearised as follows,
\begin{equation} 
	\bm{n} = \underbrace{\RMT^{*}\bm{C}_{\rm N}\bm{\Gamma}^{*}}_{\bm{n}^{\rm *}}
	- \widehat{(\RMT^{*}\bm{C}_{\rm N}\bm{\Gamma}^{*})}\Delta\RV
		+ \left(\RMT^{*}\bm{C}_{\rm N}(\RMT^{\rm *})^{\rm T}\right)
		\Delta \bm{w}^\prime + \widehat{(\bm{r}^{*})^\prime} \Delta\RV
	\label{eq:linForce}
\end{equation}
\begin{equation}
	\bm{m} = \RMT^{*}\bm{C}_{\rm M}\bm{K}^{*}
	- \widehat{(\RMT^{*}\bm{C}_{\rm M}\bm{K}^{*})} \Delta\RV 
	+ (\RMT^{*}\bm{C}_{\rm M}(\RMT^{\rm *})^{\rm T})\Delta\RV^\prime
	\label{eq:linMoment}
\end{equation}
Here, the superscript $*$ represents the fields obtained from the previous Newton-Raphson iteration and $\Delta\bm{w}$ and $\Delta\RV$ are the incremental correction displacement and rotational vectors respectively. The linearised counterpart of the term $(\bm{r}^\prime\times\bm{n})$ in Eq. (\ref{eq:momentBalance}) is given by,
\begin{equation}\label{eq:linRcrossN}
	(\bm{r}^\prime\times\bm{n}) = \widehat{(\bm{r}^*)^\prime}\bm{n}^* - 
	\widehat{\bm{n}^*}\Delta\bm{w}^\prime + \widehat{(\bm{r}^*)^\prime} \big( \bm{n} - \bm{n}^{\rm *} \big)
\end{equation}
The details of the linearisation are provided in \ref{app:linearisation}. The correction vectors $\Delta\bm{w}$ and $\Delta\RV$ are used to calculate the new mean line displacement vector $\bm{w}$ and the new rotation matrix $\RM$ at the end of each Newton-Raphson iteration according to the formulae,
\begin{equation} \label{eq:w}
	\bm{w} = \bm{w}^{*} + \Delta\bm{w}
\end{equation}
\begin{equation} \label{eq:psi}
	\RM = \exp(\widehat{\Delta\RV}) \RM^{*}
\end{equation}

\noindent where the exponentiation of the skew-symmetric tensor $\widehat{\Delta \RV} \in \mathfrak{so(3)}$  is evaluated by the Rodrigues' formula (Eq. \ref{eq:rodrigues}) to compute the rotation matrix $\RM \in SO(3)$.

\noindent \textbf{Remark}: 

\noindent The type of rotation matrix interpolation adopted in this work corresponds to the variant of Type-I as presented by Romero \cite{romero2004interpolation} and the objectivity of the corresponding strain measures is verified using a numerical test case as presented in Section \ref{subsec:objectivity}.


\section{NUMERICAL MODEL}\label{sec:numerics}

\noindent In this section, cell-centred FV discretisation of the computational domain and the governing balance equations is discussed. Unlike the standard FE approach, where the governing equation in strong form is cast into its equivalent weak form, the current FV method directly discretises the strong form of the integral equation. The discretisation procedure is separated into two distinct parts: discretisation of the solution domain and discretisation of the governing equations.

\subsection{Solution domain discretisation}\label{subsec:domDisc}
\begin{figure}[h]
	\centering
	\includegraphics[scale = 0.6]{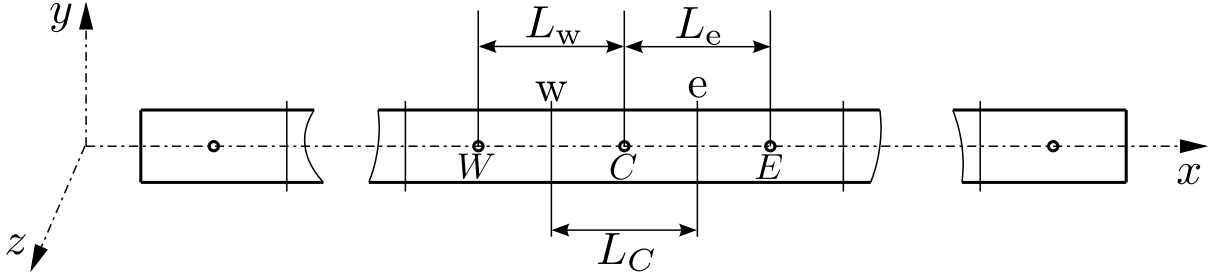}
	\caption{Beam body in reference configuration discretised by a finite set of 1-D CVs (cells).}
	\label{fig:cell}
\end{figure}

\noindent For the quasi-static case, solution domain discretisation implies space discretisation where loads are applied in pseudo-time increments. The beam body in its reference configuration is divided into a finite number of uniform segments or control volumes (CVs) as is shown in Fig. \ref{fig:cell}. A typical computational stencil (Fig. \ref{fig:cell})  consists of the central CV (cell) of length $L_{\rm C}$ with computational node $\rm C$ , located at the cell centroid, bounded by two internal faces $\rm w$ and $\rm e$ shared with the corresponding west and east neighbouring cells, with cell centroids at $\rm W$ and $\rm E$ and lengths, $L_{\rm w}$ and $L_{\rm e}$ from the node $\rm C$ respectively.

\subsection{Equation discretisation}\label{subsec:eqnDisc}

\begin{figure}[h]
	\centering
	\includegraphics[width=0.75\textwidth]{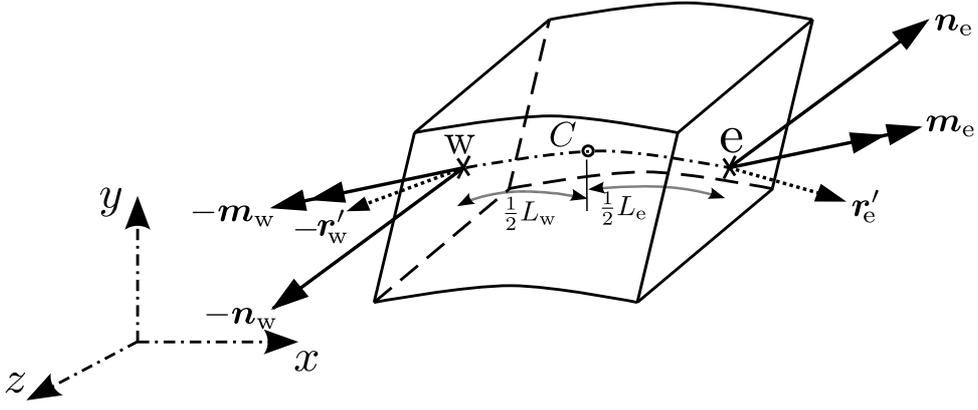}
	\caption{Balance of forces and moments on an isolated CV in the deformed configuration}
	\label{fig:CV}
\end{figure}

\noindent For an isolated CV in the deformed configuration (Fig. \ref{fig:CV}), the integral form of the balance equations (Eqs. \ref{eq:intForce} and \ref{eq:intMoment}) can be discretised over a CV as following,

\begin{align}\label{eq:discForce}
	\begin{split}
		\int_{L} \bm{n}^\prime \textrm{d}L &= \sum_{f}\bm{n}_{f} \bm{d}_{f} = \bm{n}_{\textrm{e}} - \bm{n}_{\textrm{w}} \quad ; \quad \int_{L} \bm{f} \textrm{d}L \approx \bm{f}_{\rm C} L_{\rm C} \\[0.2cm]	
		\implies &	\bm{n}_{\rm e} - \bm{n}_{\rm w} + \bm{f}_{\rm C} L_{C} = 0
	\end{split}
\end{align}

\noindent where the subscript $f$ denotes the values evaluated at cell faces and subscript $\rm C$  represents the values at cell-centre $\rm C$; $\bm{d}_{f}$ is the unit normal of the face. The term $\bm{f}$ is assumed to have a linear variation across the CV and hence, can be approximated by the mid-point rule.
 

\noindent Similarly, the discretised version of the moment balance equation about the cell-centre $\rm C$ takes the form,

\begin{align}\label{eq:discMoment}
	\begin{split}
		\int_{L} \bm{m}^\prime \textrm{d} L = \sum_{f}\bm{m}_{f} \bm{d}_f &= \bm{m}_{\textrm{e}} - \bm{m}_{\textrm{w}} \\[0.2cm]
		\int_{L} (\bm{r}^\prime \times \bm{n}) \textrm{d} L = \sum_{f} (\bm{r}^\prime_{\rm f} \times \bm{n}_{\rm f}) L_{f} \bm{d}_f  &= \frac{1}{2}L_{\rm e}(\bm{r}^\prime_{\rm e} \times \bm{n}_{\rm e})
		+ \frac{1}{2}L_{\rm w}(\bm{r}^\prime_{\rm w} \times \bm{n}_{\rm w}) \\[0.2cm]
	    \implies \bm{m}_{\rm e} - \bm{m}_{\rm w} + \frac{1}{2}L_{\rm e}(\bm{r}^\prime_{\rm e} \times \bm{n}_{\rm e})
		&+ \frac{1}{2}L_{\rm w}(\bm{r}^\prime_{\rm w} \times \bm{n}_{\rm w})
		+ \bm{t}_{\rm C} L_{C} = 0
	\end{split}
\end{align}

\noindent where, $L_{f}$ is the distance from centre $\rm C$ to faces $f$; $\bm{t}$ is also assumed to have a linear variation across the CV and hence is approximated by the mid-point rule. 


\noindent Unlike the traditional FE approach where the primary unknowns stored-at the nodal points-are used along with element shape functions to interpolate values to any other point, in the FV method, the unknown fields are computed at computational nodal points (cell-centres) and interpolated, by appropriate discretisation schemes, to the face centres and elsewhere. Subsequently, all the terms in the discretised equilibrium equations (Eqs. \ref{eq:discForce} and \ref{eq:discMoment}) are substituted by their linearised counterpart (Eqs. \ref{eq:linForce} - \ref{eq:linRcrossN}).

\vspace{0.2cm}

To evaluate the internal forces and moments at cell-faces from the linearised equations mentioned in Section \ref{sec:Math}, the mean line displacement correction vector $\Delta\bm{w}$ and its derivative $\Delta\bm{w}^\prime$, and the value of cross-section rotational correction vector $\Delta\bm{\psi}$ and its derivative $\Delta\bm{\psi}^\prime$ have to be approximated at the face centres in terms of the cell-centre values. All the cell-face derivatives are approximated by the central finite difference scheme while cell-face values are approximated by the linear interpolation. The cell-face values on the internal faces $w$ and $e$ are linearly interpolated  using cell-centre values as,

\begin{equation}\label{eq:faceValue}
	\Big[(\cdot) \Big]_{e} = \gamma_{\rm e} \Big[(\cdot) \Big]_{\rm E} + (1 - \gamma_{\rm e}) \Big[(\cdot) \Big]_{\rm C} \quad ; \quad	
	\Big[(\cdot) \Big]_{w} = \gamma_{\rm w} \Big[(\cdot) \Big]_{\rm W} + (1 - \gamma_{\rm w}) \Big[(\cdot) \Big]_{\rm C}
\end{equation}

\noindent where subscripts $\rm E$ and $\rm W$ represent values at the neighbouring cell-centres; $\Big[(\cdot) \Big]_{e}$ and $\Big[(\cdot) \Big]_{w}$ are interpolated values at internal faces $e$ and $w$ respectively and $\gamma_{\rm e}$, $\gamma_{\rm w}$	are the weighing factors given by, 

\begin{equation} \label{eq:weights}
	\gamma_{\rm e} = \frac{1}{2} \frac{L_{\rm e}}{L_{\rm C}} \quad; \quad
	\gamma_{\rm w} = \frac{1}{2}\frac{L_{\rm w}}{L_{\rm C}}
\end{equation}

\noindent The notations $\Big[(\cdot) \Big]_{\rm C}$, $\Big[(\cdot) \Big]_{\rm W}$ and $\Big[(\cdot) \Big]_{\rm E}$ denote the cell values at the location $\rm C$ and the neighbouring cells $\rm W$ and $\rm E$ (Fig. \ref{fig:cell}).  The cell face derivatives at the internal faces $\rm w$ and $\rm e$ given by  $\Big[(\cdot) \Big]^\prime_{e}$ and $\Big[(\cdot) \Big]^\prime_{w}$ are approximated using the central finite difference scheme as,

\begin{equation} \label{eq:faceDerivative}
	\Big[(\cdot) \Big]^\prime_{e} = \frac{\Big[(\cdot) \Big]_{\rm E} - \Big[(\cdot) \Big]_{\rm C} }{L_{C}} \quad ; \quad
	\Big[(\cdot) \Big]^\prime_{w} = \frac{\Big[(\cdot) \Big]_{\rm C} - \Big[(\cdot) \Big]_{\rm W} }{L_{C}}	
\end{equation}

\noindent In the discretised equilibrium equations, all the terms are \textit{explicitly} computed except for the primary unknowns, the correction vectors of incremental displacement and incremental rotation ($\Delta\bm{w}$ and $\Delta\RV$), which are treated \textit{implicitly} and are evaluated at the computational cell-centres. As mentioned in Section \ref{subsec:linearisation}, at the end of every Newton-Raphson iteration, $\Delta\bm{w}$ and $\Delta\RV$ are used to update the previously converged displacement and rotation fields (indicated by superscript $*$) and evaluate the deformed mean line position vector $\bm{r}(s)$ and the new rotation matrix $\RM$ according to the Eqs. \ref{eq:w}, \ref{eq:meanDefLine} and \ref{eq:psi}. The interpolated cell-face values and their corresponding cell-face gradients at internal faces $w$ and $e$ are used to compute the strain measures, $\bm{\Gamma}$ and $\bm{K}$ (Eqs. \ref{eq:gamma} and \ref{eq:K}), which in turn are used to evaluate the spatial forces and moments at the cell-faces (Eqs. \ref{eq:n} and \ref{eq:m}). The employed discretisation provides a nominally $2$nd order accurate approximation for displacements.

\subsection{Initial and Boundary Conditions}
\begin{figure}[h]
	\centering
	\includegraphics[width=0.6\textwidth]{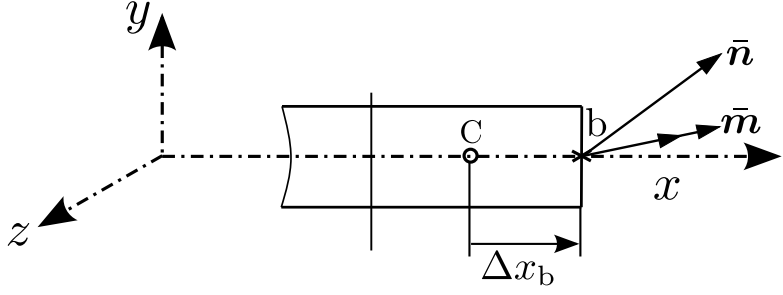}
	\caption{CV at the right boundary with a specified force $\bar{\bm{n}}$ and moment $\bar{\bm{m}}$}
	\label{fig:boundary}
\end{figure}

\noindent External forces and moments are applied in pseudo-time increments. For the beam body, there are two boundary faces, one at the left and the other at the right of the beam. 
\noindent For a Dirichlet boundary condition, in the discretised governing equations (Eqs. \ref{eq:discForce} and \ref{eq:discMoment}), the displacement/rotation component at the face centre of the boundary are directly replaced by the user defined value. For Neumann boundary conditions, the values of forces/moments at either boundary locations, are directly specified in the Eqs. \ref{eq:discForce} and \ref{eq:discMoment}. Following solution of the linear system, the corresponding displacements and rotations at that boundary location are obtained by linear extrapolation from the interiors of the solution domain using Eqs. \ref{eq:linForce} and \ref{eq:linMoment}. For instance, Fig. \ref{fig:boundary} shows the boundary face $\textrm{b}$ and the neighbouring cell-centre $\textrm{C}$. For a specified force $\bar{\bm{n}}$ and/or moment $\bar{\bm{m}}$ on the boundary, the incremental displacements/rotations can be obtained from Eqs. \ref{eq:linForce} and \ref{eq:linMoment} as,

\begin{subequations}
		\begin{align}
		\bar{\bm{n}} &= \RMT^{*}\bm{C}_{\rm N}\bm{\Gamma}^{*} + \left(\RMT^{*}\bm{C}_{\rm N}(\RMT^{\rm *})^{\rm T}\right)
		\Delta \bm{w}^\prime_{\textrm{b}} + \Big[\left(\RMT^{*}\bm{C}_{\rm N}(\RMT^{\rm *})^{\rm T}\right) \widehat{(\bm{r}^{*})^\prime} - \widehat{(\RMT^{*}\bm{C}_{\rm N}\bm{\Gamma}^{*})} \Big]\Delta\RV_{\rm b}  \\
		\bar{\bm{m}} &=\RMT^{*}\bm{C}_{\rm M}\bm{K}^{*}
		- \widehat{(\RMT^{*}\bm{C}_{\rm M}\bm{K}^{*})} \Delta\RV_{\textrm{b}} 
		+ (\RMT^{*}\bm{C}_{\rm M}(\RMT^{\rm *})^{\rm T})\Delta\RV^\prime_{\textrm{b}} \label{subeq:mom_b}
	\end{align}
\end{subequations}	 

\noindent where, the cell-face derivatives $\Delta \bm{w}^\prime_{\textrm{b}}$ and $\Delta\RV^\prime_{\textrm{b}}$ are given by,

\begin{align*}
	\Delta \bm{w}^\prime_{\textrm{b}} = \frac{\Delta \bm{w}_{\rm b} - \Delta \bm{w}_{\rm C}}{\Delta x_{\rm b}} \quad ; \quad
	\Delta \RV^\prime_{\textrm{b}} = \frac{\Delta \RV_{\rm b} - \Delta \RV_{\rm C}}{\Delta x_{\rm b}} 
\end{align*}

\noindent The value of $\Delta \RV^\prime_{\rm b}$ is substituted into Eq. \ref{subeq:mom_b} and the incremental rotation vector $\Delta \RV_{\rm b}$ is evaluated as,

\begin{equation*}
	\Delta \RV_{\rm b} = \Bigg(\frac{\RMT^{*}\bm{C}_{\rm M}(\RMT^{\rm *})^{\rm T}}{\Delta x_{\rm b}} - (\widehat{\RMT^{*}\bm{C}_{\rm M}\bm{K}^{\rm *}})\Bigg)^{-1} \Bigg[ \bar{\bm{m}} - \RMT^{*}\bm{C}_{\rm M}\bm{K}^{*} + \frac{\RMT^{*}\bm{C}_{\rm M}(\RMT^{\rm *})^{\rm T}}{\Delta x_{\rm b}} \Delta \RV_{\rm C} \Bigg]
\end{equation*}
Using $	\Delta \RV_{\rm b}$ and the value of $\Delta \bm{w}^\prime_{\textrm{b}}$ as shown above, the incremental displacement $\Delta \bm{w}_{\rm b}$ is then calculated as,

\begin{equation*}
	\Delta \bm{w}_{\rm b} = \Delta \bm{w}_{\rm C} + \Bigg(\frac{\RMT^{*}\bm{C}_{\rm N}(\RMT^{\rm *})^{\rm T}}{\Delta x_{\rm b}}\Bigg)^{-1} \Bigg[ \bar{\bm{n}} - \RMT^{*}\bm{C}_{\rm N}\bm{\Gamma}^{*} - \Big[\left(\RMT^{*}\bm{C}_{\rm N}(\RMT^{\rm *})^{\rm T}\right) \widehat{(\bm{r}^{*})^\prime} - \widehat{(\RMT^{*}\bm{C}_{\rm N}\bm{\Gamma}^{*})} \Big]\Delta\RV_{\rm b} \Bigg]
\end{equation*}

\subsection{Solution procedure}

\noindent The final form of the discretised equilibrium equations, with appropriate discretisation schemes described in Section \ref{subsec:eqnDisc}, for a typical computational node $\rm C$ read as follows,
\begin{equation}
	\bm{A}_{\rm W}	
	\begin{bmatrix}
		(\Delta\bm{w})_{\rm W}\\
		(\Delta\bm{\psi})_{\rm W}
	\end{bmatrix}
	+ \bm{A}_{\rm C}	
	\begin{bmatrix}
		(\Delta\bm{w})_{\rm C}\\
		(\Delta\bm{\psi})_{\rm C}
	\end{bmatrix}
	+ \bm{A}_{\rm E}	
	\begin{bmatrix}
		(\Delta\bm{w})_{\rm E}\\
		(\Delta\bm{\psi})_{\rm E}
	\end{bmatrix}
	= \begin{bmatrix}
		(\bm{R}^{\bm{w}})_{\rm C}\\
		(\bm{R}^{\RV})_{\rm C}
	\end{bmatrix} 	
	\label{eq:discreteEquation}
\end{equation}

\noindent where $\bm{A}_{\rm C}$ is a coefficient matrix containing the contributions of node $\rm C$ while the matrices $\bm{A}_{\rm W}$ and $\bm{A}_{\rm E}$ representing the interactions with the neighbouring cell centres $\rm W$ and $\rm E$. The right hand side of the Eq. \ref{eq:discreteEquation} is the source vector contribution. All the coefficient matrices are $(6 \times 6)$ dense coupled matrices with the primary unknowns being $\Delta \bm{w}$ and $\Delta \RV$. The three components of the mean line displacement correction and cross-section rotation vectors have to be solved in coupled manner. The detailed structure of the diagonal and off-diagonal coefficient matrices is provided in the \ref{app:BlockCoeffs}. 

The linearised Eqs. \ref{eq:discreteEquation} are assembled for all CVs forming a system of equations given by,

\begin{equation}\label{eq:assemble}
	\big[ \bm{A} \big] \big[\bm{\phi}\big] = \big[\bm{R}\big]
\end{equation}

\noindent resulting in $6M \times 6M$ sparse matrix $[\bm{A}]$ with weak diagonal dominance, with $M$ being the total number of CVs. The coefficients  $\bm{A}_{\rm C}$ constitute the diagonal whereas matrices $\bm{A}_{\rm W}$ and $\bm{A}_{\rm E}$ contribute to the off-diagonal terms of $[\bm{A}]$. The solution vector $\big[\bm{\phi}\big]$ contains the primary unknowns $\Delta \bm{w}$ and $\Delta \RV$, and $\big[\bm{R}\big]$ is the source vector containing the explicit discretised terms and boundary condition contributions. The final system of linearised algebraic equations obtained by assembling Eqs. \ref{eq:discreteEquation}  for all control volumes in the mesh is successfully solved using the block variant of the Thomas algorithm.

For every pseudo-time increment, the coupled equations are iteratively solved by Newton-Raphson procedure, until a convergence tolerance of order $10^{-10}$ is achieved with the maximum number of allowable iterations set to $30$. For calculating the normalised (with respect to their magnitudes) residuals, after each iteration, for all CVs, a list of the magnitudes of the normalised correction vectors, $\Delta \bm{w}$ and $\Delta \RV$ is calculated and the global maximum of the two magnitudes is selected ($\textrm{Residuals} = max\Big[ |\Delta \bm{w}|, |\Delta \RV| \Big]$) and the equations are solved until the convergence tolerance is achieved. The current method has been implemented in open-source software OpenFOAM (\cite{weller1998tensorial}) (version foam-extend-4.1), exploiting the developed object oriented FV procedures. The overall solution procedure is summarised in Algorithm \ref{algo:solution} and the procedure to update the kinematic quantities and stress resultants is given in Algorithm \ref{algo:update}.

\begin{algorithm} 
	\caption{Solution procedure; \qquad \qquad \qquad \qquad \qquad \qquad \qquad \qquad \qquad $(\cdot)_{f}$ : fields at the face centres}
	\begin{algorithmic}[1] 
		\State Set 
		 $\RMI = \mathbf{I}$  	\Comment{Reference mean tangent line and reference rotation matrix} 
		\If{$\textrm{initially curved beam}$}
		\State Set $\bm{r}_0$ and $\RMI$ at cell faces 		 \Comment{Initial undeformed mean line and initial rotation matrix; specific to test case} 
		\Statex \qquad $\implies (\bm{r}_0)_f, (\RMI)_f$
		\EndIf
		\Statex \textbf{end if}
		\State Set $\RM = \mathbf{I}$, $\RM_f = \mathbf{I}$ and $\RMT = \RM (\RMI)_f$ 	\Comment{Relative rotation and total rotation matrix}
		\vspace{0.1cm}
		\State Calculate $\big(\bm{r}_0^\prime\big)_f$ from $(\bm{r}_0)_f$
		\vspace{0.1cm}
		\State Set $\bm{g}_0 =  \big(\bm{r}_0^\prime\big)_f$ 							\Comment{Initial tangent vector at face centres}
		\For{all pseudo-time steps}
		\While{residuals are not converged}
		\State Coupled linearised balance equations: assemble and solve for $\Delta \bm{w}$ and $\Delta \bm{\psi}$ \Comment{Eq. \ref{eq:assemble}}
		\State Interpolate cell-centre displacements $\Delta \bm{w}$ and rotation vectors $\Delta \bm{\psi}$  to face values \Comment{Eq. \ref{eq:faceValue}}
		\State Update kinematic quantities and stress resultants  \Comment{Algorithm \ref{algo:update}}
		\EndWhile 
		\State \textbf{end while} 
		\EndFor
		\State \textbf{end for} 
	\end{algorithmic}
	\label{algo:solution}
\end{algorithm}

\begin{algorithm} 
	\caption{Update of the kinematic quantities and stress resultants;  \quad \quad  $(\cdot)_{f}$ : fields at the face centres \\ $(\cdot)^*$ : fields calculated in the previous Newton-Raphson iteration}
	\begin{algorithmic}[1]
		\State $ \bm{w} = \bm{w}^* + \Delta \bm{w}$				\Comment{Update displacement vector (Eq. \ref{eq:w})}
		
		\State $\RV = \RV^* + (\RM^{*})^{T} \Delta \RV$				\Comment{Update rotation vector}
		
		\State Calculate $(\Delta \RV)_f$ from $\Delta \RV$ 	\Comment{Interpolate rotation vector to the face centre (Eq. \ref{eq:faceValue})}
		
		\vspace{0.1cm}
		
		\State $\Delta \bm{T}_f = \frac{\sin \|\Delta \RV_f\| }{\|\Delta \RV_f\|} \mathbf{I} + \Big(1 - \frac{\sin \|\Delta \RV_f\|}{\|\Delta \RV_f\|}\Big) \frac{\Delta \RV_f \cdot \Delta \RV^T_f }{\|\Delta \RV_f\|^2} + \Big( \frac{1 - \cos \|\Delta \RV_f\|}{\|\Delta \RV_f\|^2} \Big) \widehat{\Delta \RV_f}^2$				\Comment{Calculate incremental tangent operator (Eq. \ref{eq:tangent})}
		
		\Statex
		
		\State $\bm{K}_f = \bm{K}_f^* + \big(\RMI^{T} \big)_f (\RM^{*})^{T}_f \Delta \bm{T}_f^T (\Delta \RV^\prime)$ \Comment{Calculate (material) rotational strains (Eq. \ref{eq:K})}
		
		
		\State $(\Delta \RM)_f = \mathbf{I} + \Big(\frac{\sin \|\Delta \RV_f\|}{\|\Delta \RV_f\|}\Big) \Delta \RV_f \cdot \Delta \RV^T_f  + \Big( \frac{1 - \cos \|\Delta \RV_f\|}{\|\Delta \RV_f\|^2} \Big) \widehat{\Delta \RV_f}^2$ \Comment{Calculate face centre incremental rotation matrix (Eq. \ref{eq:rodrigues})}
		
		\Statex
		
		\State $\RM_f = (\Delta \RM)_f \RM_f^*$ \Comment{Update the face centre rotation matrix (Eq. \ref{eq:psi})}
		
		\vspace{0.1cm}
		
		\State $\Delta \RM = \mathbf{I} + \Big(\frac{\sin \|\Delta \RV\|}{\|\Delta \RV\|}\Big) \Delta \RV \cdot \Delta \RV^T  + \Big( \frac{1 - \cos \|\Delta \RV\|}{\|\Delta \RV\|^2} \Big) \widehat{\Delta \RV}^2$ \Comment{Calculate the incremental cell centre rotation matrix (Eq. \ref{eq:rodrigues})}
		\vspace{0.1cm}	
		
		\State $\RM = (\Delta\RM) \RM^*$	\Comment{Update the cell centre rotation matrix}
		
		\vspace{0.1cm}
		
		\State $\big(\bm{r}^\prime \big)_f = \big(\bm{r}_0 \big)_f + \bm{w}^\prime $ \Comment{Calculate  derivative of deformed mean line at face centre}
		
		\vspace{0.1cm}
		
		\State $\bm{\Gamma}_f = \big(\RMI^T\big)_f \RM_f^T \big(\bm{r}^\prime \big)_f - \big(\RMI^T\big)_f \big(\bm{r}^\prime_0 \big)_f $
		\vspace{0.1cm}	\Comment{Calculate (material) translational strains (Eq. \ref{eq:gamma})}
		
		\State $\bm{n}_f = \RM_f \big(\RMI \big)_f \bm{C}_{\rm N} \bm{\Gamma}$ \Comment{Calculate spatial forces at face centres (Eq. \ref{eq:n})}
		
		\vspace{0.1cm}
		
		\State $\bm{m}_f = \RM_f \big(\RMI \big)_f \bm{C}_{\rm M} \bm{K}$	\Comment{Calculate spatial moments at face centres (Eq. \ref{eq:m})}
 
		
	\end{algorithmic}
	\label{algo:update}
\end{algorithm}


\section{VERIFICATION TEST CASES}\label{sec:testCases}

\noindent In this section, the capabilities of the developed FV methodology are investigated on the five complementary benchmark cases:
\begin{enumerate}
	\item Rigid rotation of an initially curved beam: this case verifies the objectivity of the adopted strains and rotation interpolation measures.
	\item In-plane bending of a cantilever beam subjected to concentrated moment at one end: this case provides an analytical solution to check the accuracy  and order of accuracy of the proposed methodology. 
	\item Out-of-plane bending of the cantilever to form a helix due to a concentrated moment and an out-of-plane force: this case with complex loading conditions, tests the ability of the numerical solver to cater for large rotations and 3-D deformation.
	\item A cantilever, initially bent into a $45^{\circ}$ arc in $xy$-plane, subjected to a force along the $z$-direction: this is another benchmark curved beam case used by many authors to establish the accuracy and order of accuracy of the numerical model.
	\item Deep-circular arch with a concentrated in-plane force at the crown location: this unsymmetrical circular arch case has exact solutions available and is used to assess the ability of the solver to obtain the critical buckling load value.
\end{enumerate}
\noindent All the test cases have been executed using a quad core CPU with hyper-threading (Intel Core(TM) i7-8565U CPU with base frequency $1.80$GHz and maximum turbo frequency $4.6$GHz).

\subsection{Rigid rotation of an initially curved beam} \label{subsec:objectivity}

\begin{figure}[h]
	\centering
	\begin{subfigure}{0.4\textwidth}
		\centering
		\includegraphics[scale = 0.5]{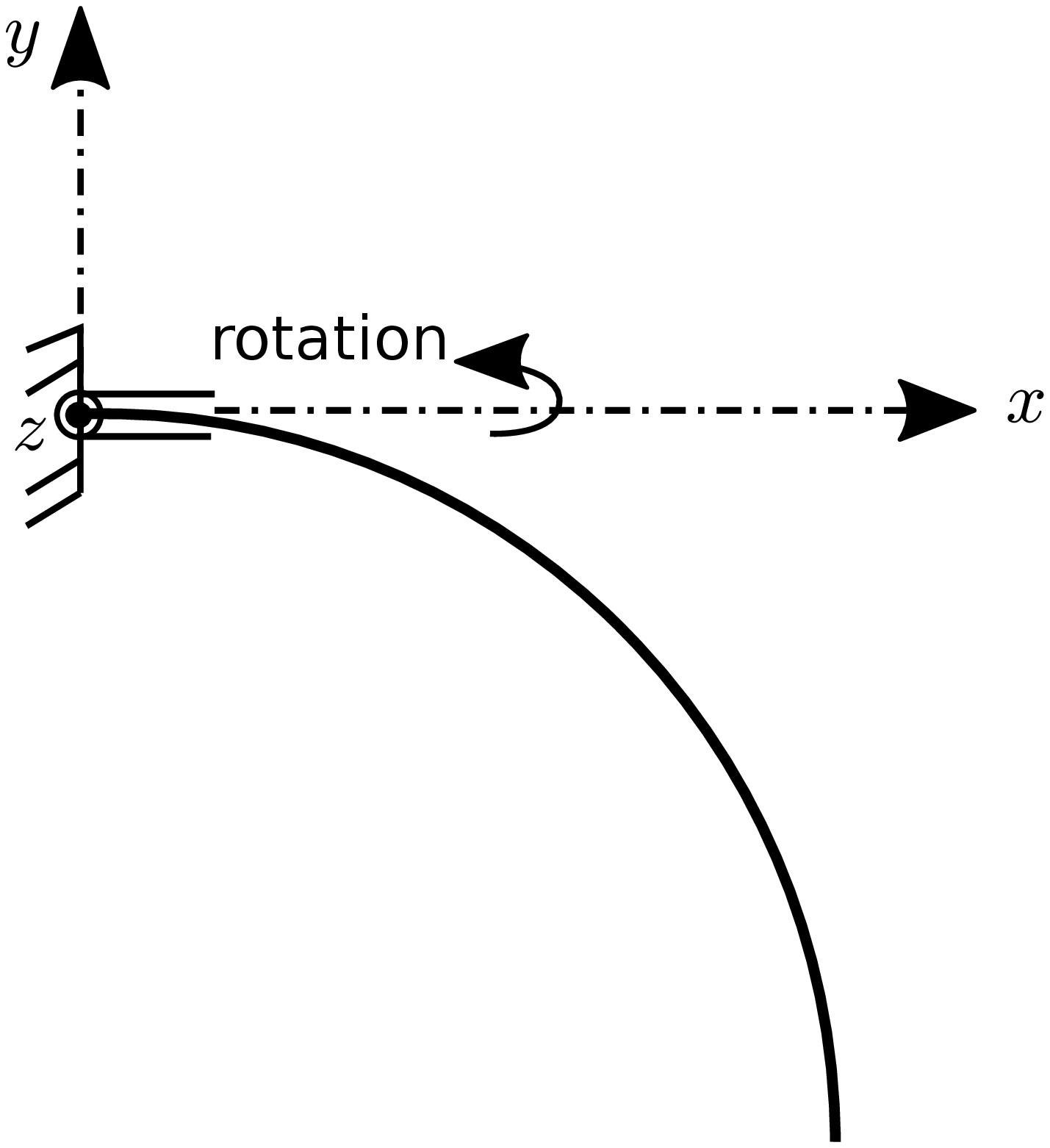}
		\caption{}
		\label{fig:objective}
	\end{subfigure}
	\hfill
	\begin{subfigure}{0.55\textwidth}
		\centering
		\includegraphics[scale = 0.5]{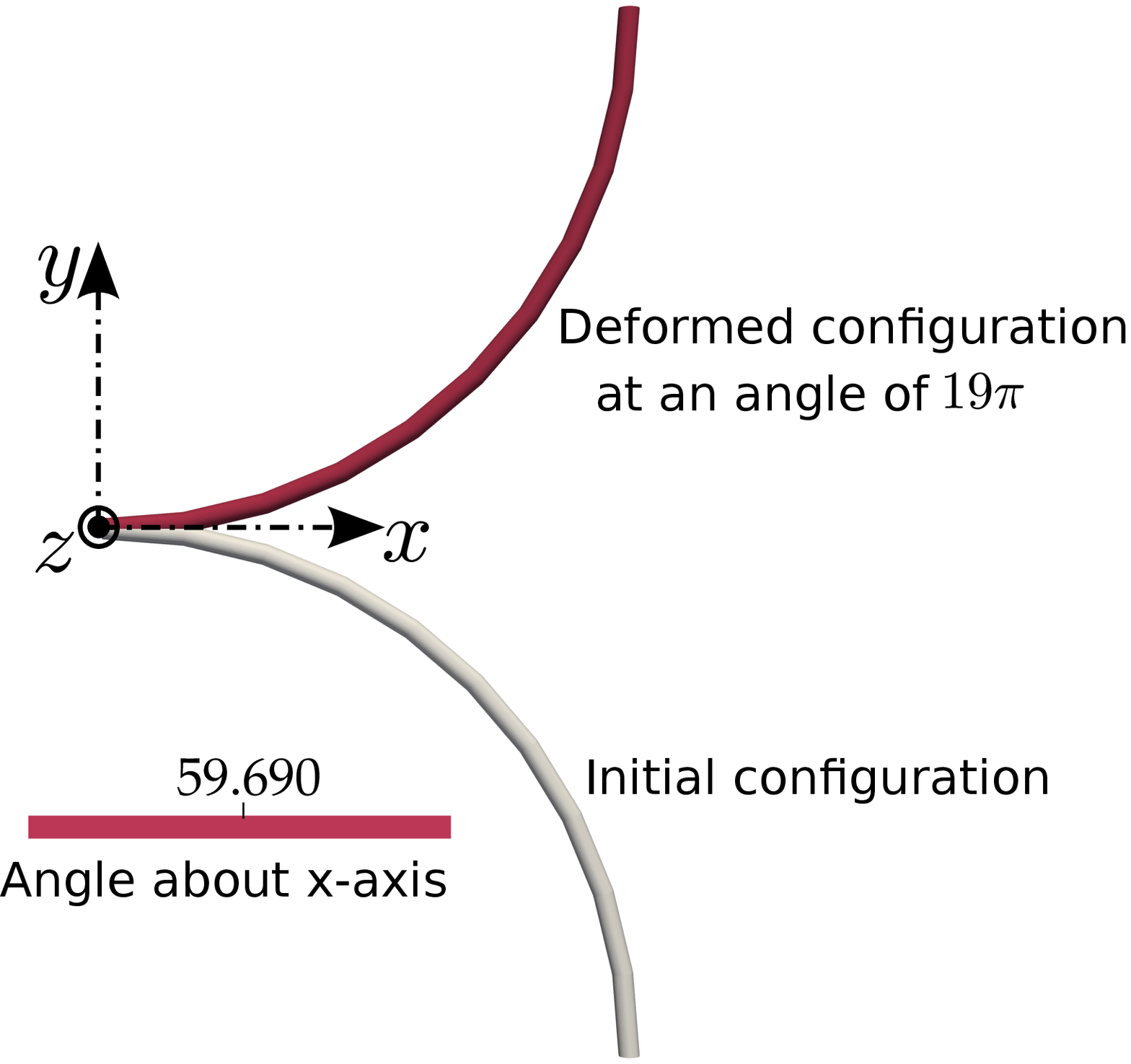}
		\caption{}
		\label{fig:objectiveresult}
	\end{subfigure}
	\caption{Rigid rotation of an initially curved beam: (a) test case setup, (b) deformed configuration of the beam for a rotation angle, $\psi_x = 19\pi$}
	\label{fig:objectiveTestcase}
\end{figure}

\noindent For the FV formulation presented in this paper and the type of the rotation interpolation adopted, the objectivity of the conjugate strain measures is verified using the numerical test case first presented in Meier \textit{et al.} \cite{meier2014objective,meier2019geometrically} as shown in Fig. \ref{fig:objective}. For an initially curved beam with a centerline configuration of a quarter circle with radius $R = 100$ m, initial length, $L = \pi R/2 = 157.079$ m, discretised into $10$ CVs, a rotation of $20\pi$ is gradually applied about the global x-axis in $100$ load increments. The mechanical properties for this test case are, $E = 1$ GPa and $G = 0.5$ GPa. The nature of the boundary condition at the clamped left end of the beam is such that, this prescribed rotation should only cause rigid body transformation in the beam about the $x$-axis with no accumulation of internal strain energy and no physical deformation, i.e. $w_x = w_y = w_z = 0$ and $\psi_y = \psi_z = 0$.

Fig. \ref{fig:objectiveresult} shows the deformed beam configuration for an angle $\psi_x=19\pi$ about the x-axis and it is evident that there is no physical deformation. For the test case, execution time is approximately $1.5$ s and requires $9$ outer iterations on average to converge per pseudo-time increment. For each CV, the strains $\bm{\Gamma}$ and $\bm{K}$ are used to calculate the internal energy as per Eq. \ref{eq:intEnergy} and summed over all the CVs to obtain the total  accumulated strain energy of the beam due to the applied rotation. The total energy is found to be zero to machine precision (in the order of $10^{-31}$) and hence, the objectivity of the adopted strain measures is confirmed.

\subsection{In-plane bending of a cantilever beam subjected to concentrated moment at one end}

\noindent This in-plane pure flexural bending of a cantilever case has been investigated by a number of authors~\cite{simo1986three,ibrahimbegovic1995computational}. An initially straight cantilever of length $L = 10$ $\textrm{m}$ is bent into a circle by applying a concentrated moment at one end. The mechanical properties available in the literature for this test case are $EA = 10^4$ $\textrm{N}$ , $ GA_{2} = GA_{3} = 5000$ N, $EI_{2} = EI_{3} = 100$ $\textrm{Nm}^2$, $GJ = 100$ $\textrm{Nm}^2$ respectively. Hence, the cross-section radius of the beam, $r = 0.2$ m and the Young's modulus $E = 7.95 \times 10^{4}$ Pa are assumed in a way to achieve these desired numerical values. The Poisson's ratio $\nu$ is taken as zero. According to the classic Euler formula, the analytical solution for pure beam-bending is given by,

\begin{equation}\label{eq:circleAnalytic}
	\psi_z = \frac{M_z L}{EI} \quad;\quad w_{x} = L - \frac{L}{\psi_z/2} sin\frac{\psi_z}{2} cos\frac{\psi_z}{2} \quad;\quad w_{y} = \frac{L}{\psi_z/2} \Big( sin\frac{\psi_z}{2} \Big)^2  
\end{equation}
where $\psi_z$, $w_{x}$ and $w_{y}$ are the rotation and in-plane displacements of the beam tip. 

\begin{figure}[h]
	\centering
	\includegraphics[width=0.6\textwidth]{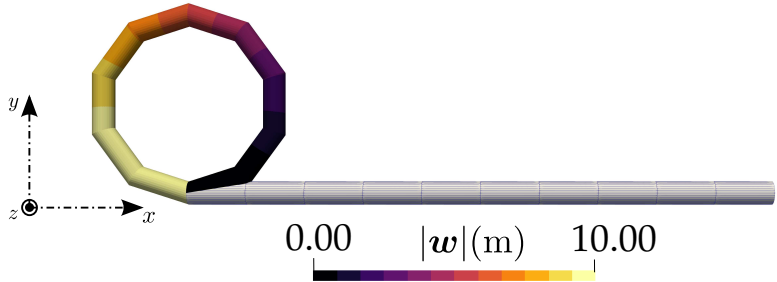}
	\caption{In-plane bending of a cantilever beam subjected to concentrated moment at one end: deformed configuration for $\psi_z = 2\pi$}
	\label{fig:cantileverCircle}
\end{figure}

\begin{figure}[h]
	\centering
	\includegraphics[width=0.8\textwidth]{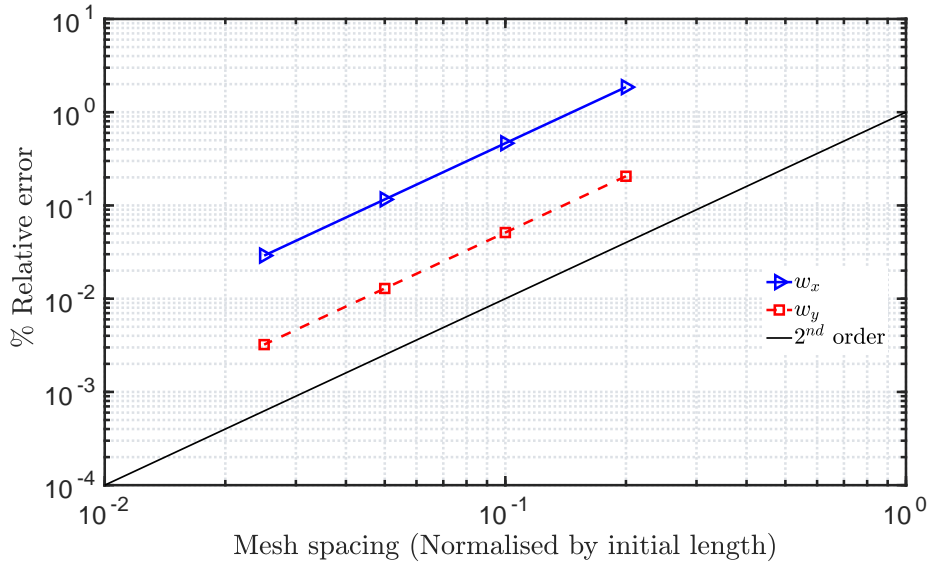}
	\caption{In-plane bending of a cantilever beam subjected to concentrated moment at one end: mesh convergence plot of in-plane displacements for an applied moment $M = 2.5\pi$ $\textrm{Nm}$}
	\label{fig:cantileverCircleConv}
\end{figure}

The computational domain is discretised into $10$ CVs and a moment of $M_z = 20\pi$ Nm $(\psi_z = 2\pi)$ is applied at the right end in $4$ increments and the Newton-Raphson solution procedure converges in an average of $6$ outer iterations per pseudo-time increment. Fig. \ref{fig:cantileverCircle} shows the deformed configuration of the cantilever for the applied moment. Although the solver takes larger number of iterations to converge than is reported in the literature (Simo et al. \cite{simo1986three} reported only two iterations), the execution time is less than $0.5$ s. For an applied moment, $M_z = 2.5\pi$ Nm and a mesh of $10$ CVs, the tip in-plane displacements are found to be $w_{x} = -1.00146$ m and $w_y = 3.72731$ m respectively, which vary from the analytical solution (Eq. \ref{eq:circleAnalytic}) by only $0.4\%$ and $0.05\%$ respectively. The percentage relative error is calculated by the formula,

\begin{equation}\label{eq:relError}
	\textrm{\% relative error} = \bigg| \frac{\xi^{\textrm{num}} - \xi^{\textrm{ref}}}{\xi^{\textrm{ref}}} \bigg| \times 100\%
\end{equation}

\noindent where, $\xi^{\textrm{num}}$ denotes the numerical value and $\xi^{\textrm{ref}}$ is the reference/analytical result. Fig. \ref{fig:cantileverCircleConv} shows the percentage mesh error convergence of the in-plane displacements for successive reduction of mesh sizes, i.e. $5, 10, 20$ and $40$ CVs; a quadratic order of error convergence is observed.

\subsection{Out-of-plane bending of the cantilever to form a helix due to a concentrated moment and an out-of-plane force}

\noindent For this case, the previous problem is extended by applying a concentrated out-of-plane force at the free end of the beam, $\bm{n} =[0, 0, 50]$ N along with a rotation $\psi_z = 20\pi$  ($\bm{m} =[0, 0, 200\pi]$ Nm) about the z-axis. Thus, an initially straight cantilever beam is bent to a circular helix shape. The constitutive matrices are taken similar to the previous case, i.e. $\bm{C}_{\rm N} = \diag[10^4, 5 \times 10^3, 5 \times 10^3]$ N and $\bm{C}_{\rm M} = \diag[10^2, 10^2, 10^2]$ Nm$^2$. The cross-section radius is reduced to $r = 0.02$ m and the Young's modulus is taken as, $E = 7.95 \times 10^{8}$ Pa, so that the beam does not overlap into its own body while forming the helix. For comparing the achieved numerical results with those published by Ibrahimbegovic \cite{Ibrahimbegovic199749}, the beam of initial length $L = 10$ m is discretised into $100$ uniform CVs. 
 
 \begin{figure}[t]
 	\centering
 	\begin{subfigure}{\textwidth}
 		\centering
 		\includegraphics[width=\textwidth]{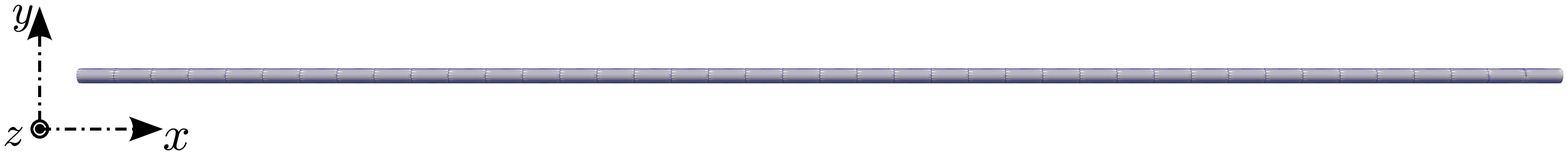}
 		\caption{}
 		\label{fig:initialBeam}
 	\end{subfigure}
 	\hfill
 	\begin{subfigure}{\textwidth}
 		\centering
 		\includegraphics[width = \textwidth]{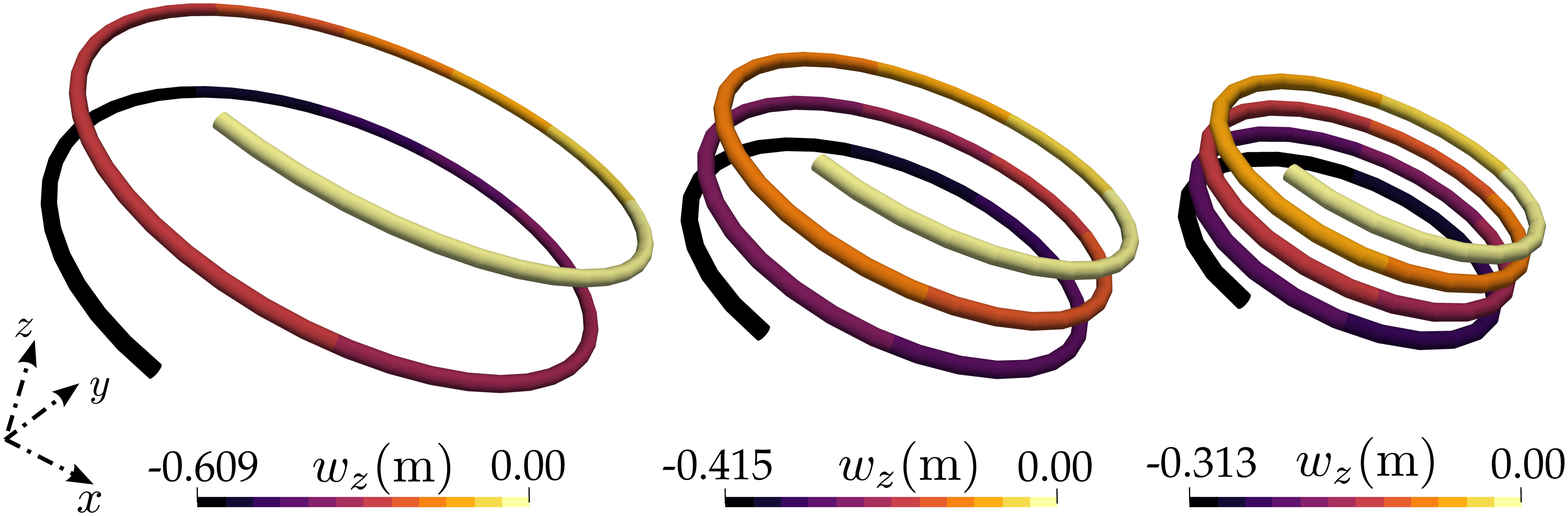}
 		\caption{}
 		\label{fig:deformedBeam}
 	\end{subfigure}
 	\caption{Out-of-plane bending of the cantilever to form a helix due to a concentrated moment and an out-of-plane force: (a) initial beam configuration (b) deformed beam for $20\%$ (left), $30\%$ (middle) and $40\%$ (right) of the applied load}
 	\label{fig:helix}
 \end{figure}

\begin{figure}[!]
	\centering
	\includegraphics[width=0.8\textwidth]{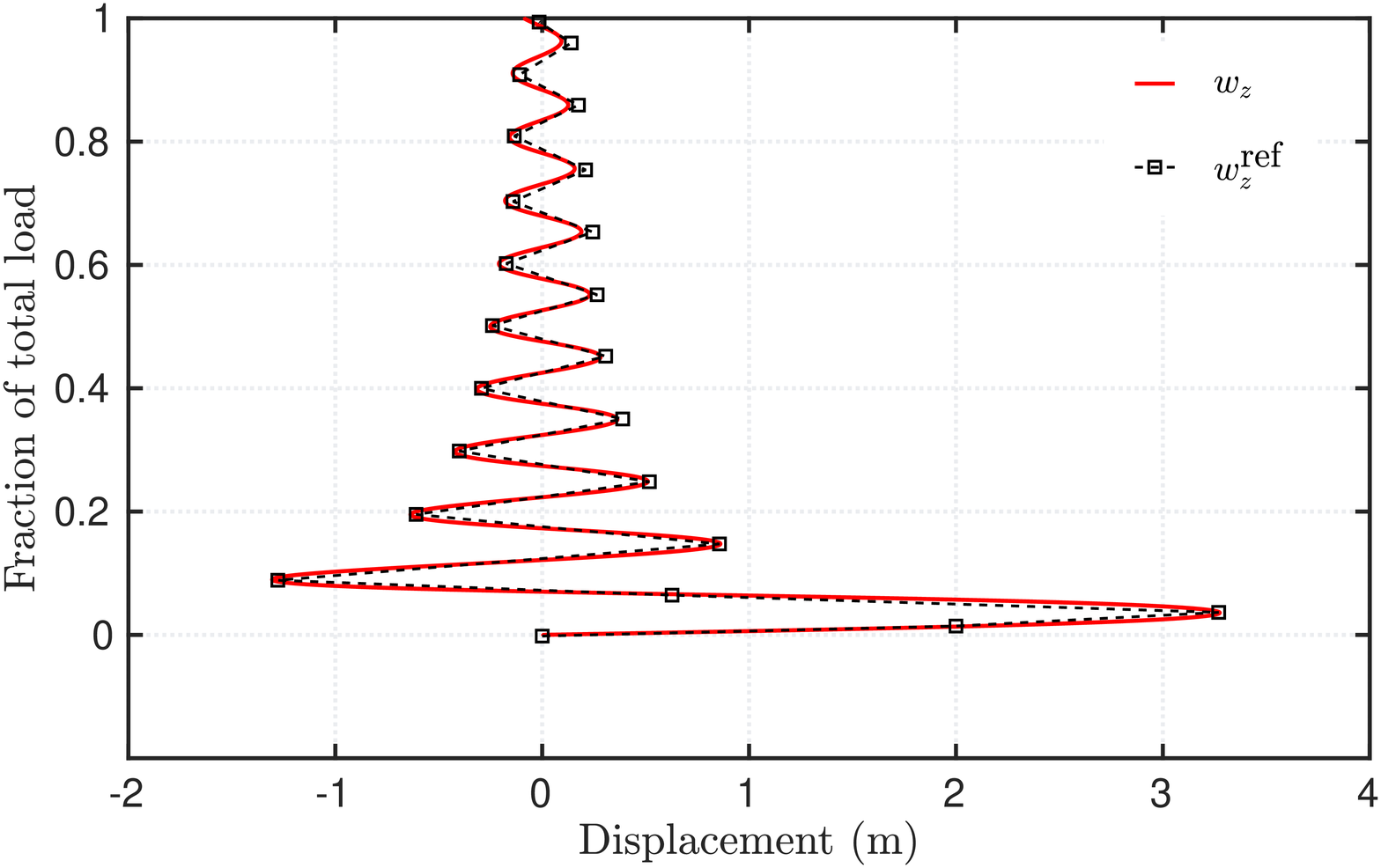}
	\caption{Out-of-plane bending of the cantilever to form a helix due to a concentrated moment and an out-of-plane force: free-end displacement component in the direction of applied load ($w^{\textrm{ref}}_z$ adopted from \cite{Ibrahimbegovic199749})}
	\label{fig:helixRes}
\end{figure}

\noindent The external loads to the beam are applied in $2000$ pseudo-time increments and the final load is reached at pseudo time-value of $10$ s. The Newton-Raphson solution procedure required $6$ outer iterations on average to converge per pseudo time increment and $60$ s of total execution time. Fig. \ref{fig:initialBeam} shows the initial beam and Fig. \ref{fig:deformedBeam} shows the deformed shape of the beam for $20\%$, $30\%$ and $40\%$ of the total loading respectively. Fig. \ref{fig:helixRes} presents the free-end displacement component $w_z$ as a function of the applied load. From Fig. \ref{fig:helixRes}, it is evident that with increasing rotations, the out-of-plane displacement oscillates about the $z$-axis crossing the zero value. This ability of the solver to capture the complex deformation and oscillatory motion of the beam is possible because of the incremental rotation vector based parametrisation and the results (Fig. \ref{fig:helixRes}) are in good agreement with those published in \cite{Ibrahimbegovic199749}; note the same test case failed to converge when a non-incremental formulation was initially trialed.

\subsection{A cantilever, initially bent into a $45^{\circ}$ arc in $xy$-plane, subjected to a force along the $z$-direction}

\begin{figure}[ht]
	\centering
	\includegraphics[scale = 0.4]{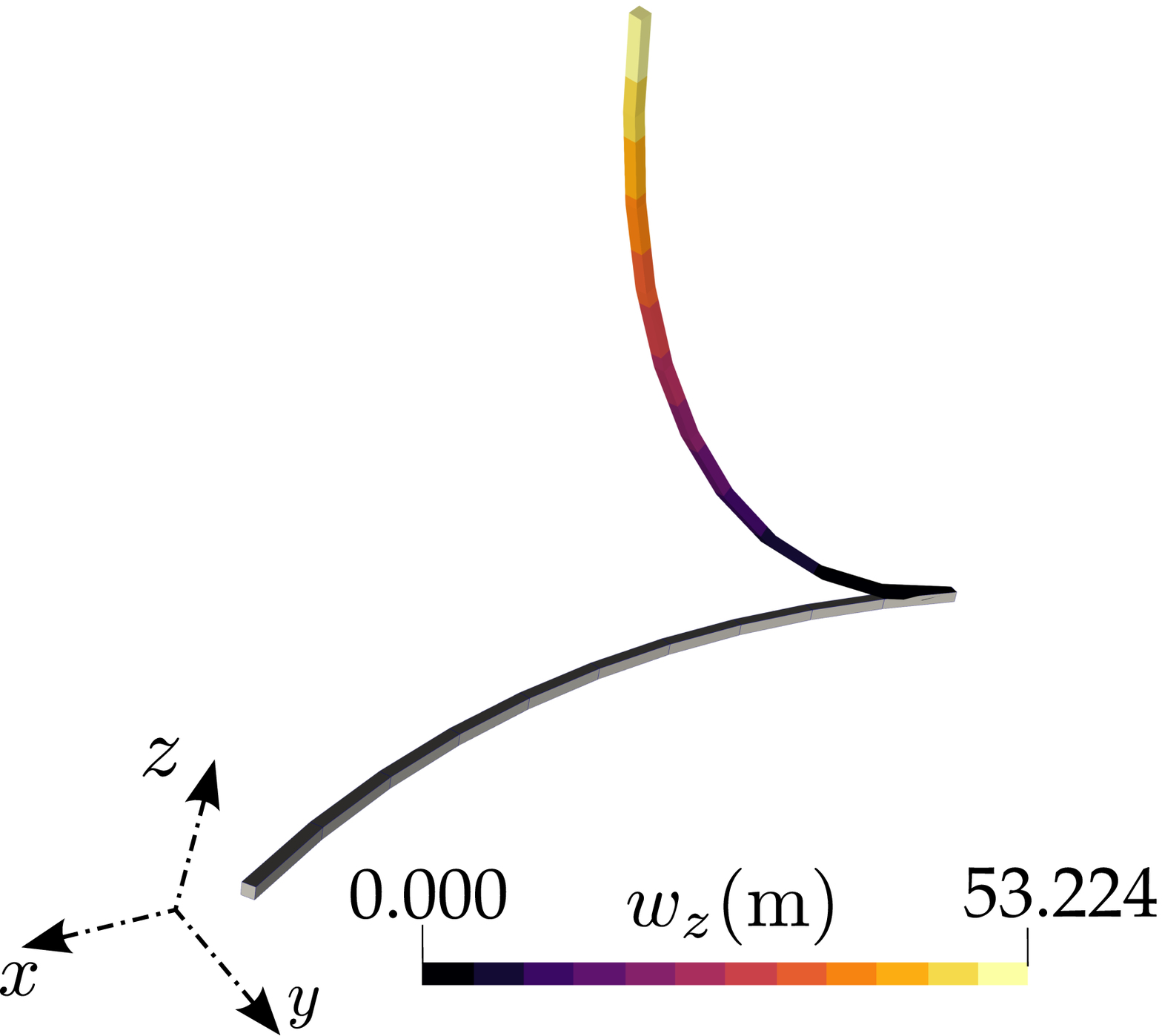}
	\caption{A cantilever, initially bent into a $45^{\circ}$ arc in $xy$-plane, subjected to a force along the $z$-direction: initial and deformed configuration for an applied load of $600$ N}
	\label{fig:45degArc}
\end{figure}

\noindent This case consists of an initially curved cantilever beam with an applied vertical concentrated force along $z$-axis. Several authors have reproduced this benchmark case \cite{bathe1979large,simo1986three,cardona1988beam,crisfield1990consistent,ibrahimbegovic1995computational}. The cantilever, which has a unit square cross-section, is initially bent into a $45^{\circ}$ arc of radius $100$ m, length $L = \pi R/4 = 78.54$ m and then a vertical force, $n_{z} = 600$ N is applied. The mechanical properties adopted are $E = 1 \times 10^7$ $Pa$ and $\nu = 0$. For a spatial discretisation of $8$ CVs, the tip end displacements values are $w_x = -23.5281$ m, $w_y = -13.5415$ m and $w_z = 53.082$ m respectively. Simo et al. \cite{simo1986three} apply a load in three steps ($300$, $150$ and $150$  N) and take 27 cumulative iterations to attain convergence. The current model takes $34$ cumulative iterations for the same loading increment but less than $0.5$ s of execution time. On the other hand, for a systematic loading in $10$ increments and discretisation of $10$ CVs, it take an average of $6$ outer iterations to converge per pseudo time increment. Fig. \ref{fig:45degArc} shows the initial and the deformed configuration of the beam and the displacement component $w_z$ for $10$ CVs.

\begin{table}[ht]
	\centering
	\caption{Displacement values observed for a curved cantilever in $xy-$plane subjected to a force along $z-$axis}
	\label{table:45DegArcResults}       
	\begin{tabular}{cccc}
		\hline\noalign{\smallskip}
		Numerical Results & $|w_{x}|$ & $|w_{y}|$  & $|w_{z}|$\\
		\noalign{\smallskip}\hline\noalign{\smallskip}
		Present	& $23.540$ & $13.564$ & $53.225$\\
		Bathe and Bolourchi\textsuperscript{ \cite{bathe1979large}} 	& $23.5$ & $13.4$ & $53.4$\\
		Simo and Vu-Quoc\textsuperscript{ \cite{simo1986three}} 	& $23.48$ & $13.50$ & $53.37$\\
		Cardona and Geradin\textsuperscript{ \cite{cardona1988beam}} 	& $23.67$ & $13.73$ & $53.50$\\
		Crisfield\textsuperscript{ \cite{crisfield1990consistent}} 	& $23.87$ & $13.63$ & $53.71$\\
		Ibrahimbegovi\'{c}\textsuperscript{ \cite{ibrahimbegovic1995computational}} 	& $23.697$ & $13.668$ & $53.498$\\[0.5ex]
		\noalign{\smallskip}\hline
	\end{tabular}
\end{table}

\begin{figure}[!]
	\centering
	\includegraphics[width=0.8\textwidth]{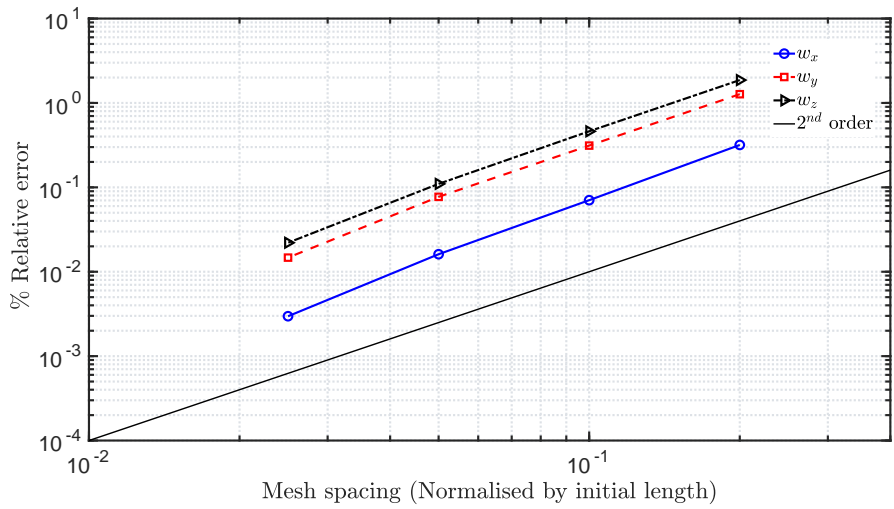}
	\caption{A cantilever, initially bent into a $45^{\circ}$ arc in $xy$-plane, subjected to a force along the $z$-direction: mesh convergence of the displacement components for the applied load, $n_z = 600$ N}
	\label{fig:45DegArcMeshConv}
\end{figure}

\begin{figure}[!]
	\centering
	\includegraphics[width=0.8\textwidth]{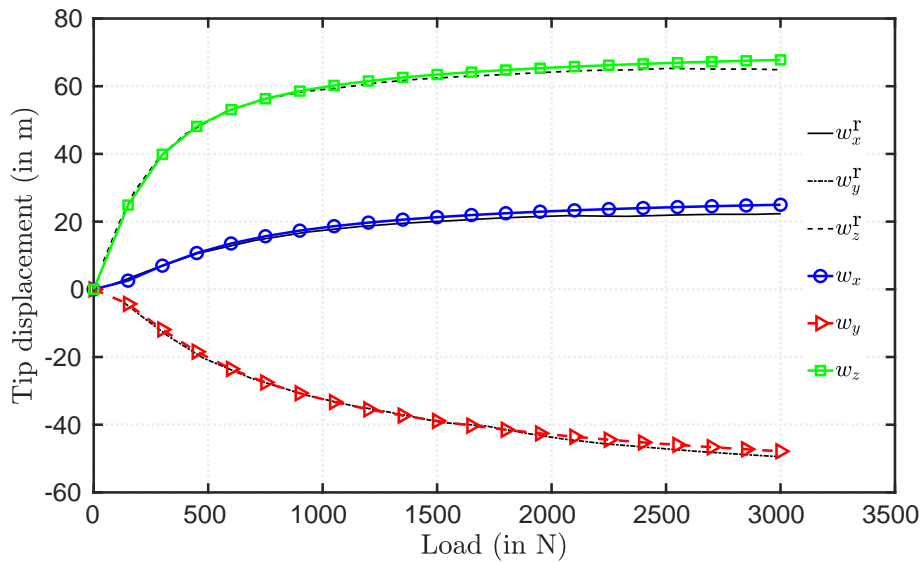}
	\caption{A cantilever, initially bent into a $45^{\circ}$ arc in $xy$-plane, subjected to a force of $3000$ N along the $z$-direction: the end-tip displacements for increasing load values}
	\label{fig:45degArcLoad}
\end{figure}

\begin{figure}[!]
	\centering
	\includegraphics[width=0.8\textwidth]{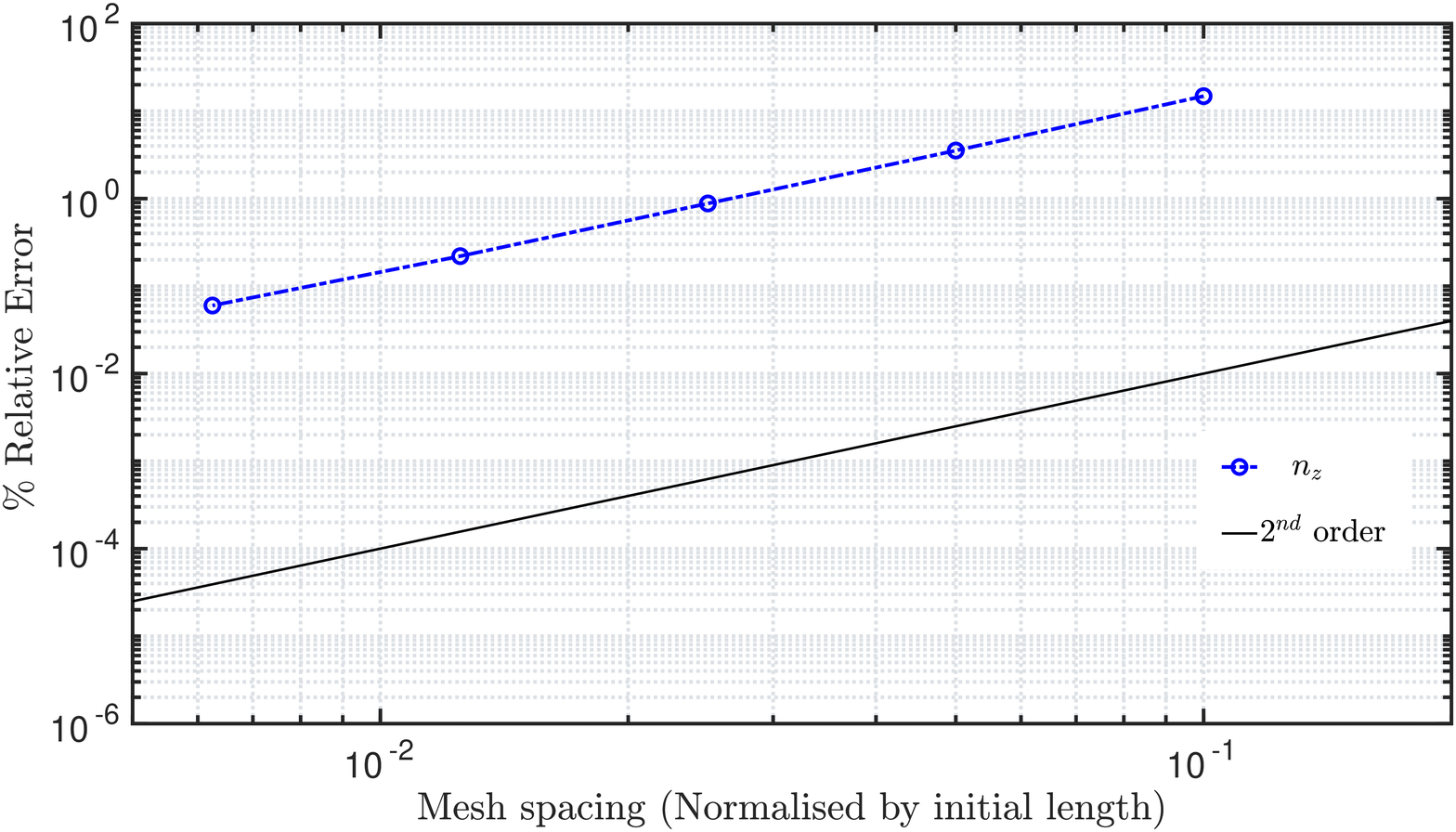}
	\caption{A cantilever, initially bent into a $45^{\circ}$ arc in $xy$-plane: mesh convergence of the force $n_z$ for an applied displacement, $w_x = -23.5607$ m, $w_y = -13.6048$ m, $w_z = 53.4756$ m}
	\label{fig:45DegArcMeshConvForce}
\end{figure}

\noindent For this test case, the reported displacements in the literature are for a discretisation of $8$ linear FE elements and the current numerical results (for $10$ CVs) are found to be comparable with the ones presented in the literature (Table \ref{table:45DegArcResults}). Fig. \ref{fig:45DegArcMeshConv} presents the \% mesh error convergence of the displacement components for successive mesh size reductions, viz. $5, 10, 20$ and $40$ CVs; a quadratic order of error convergence is observed. The analytical solutions for this test case are not available and the numerical results reported in literature are similar. To calculate the \% mesh error (Eq. \ref{eq:relError}), the displacement values corresponding to a finer mesh of $80$ CVs are adopted as reference results here. For comparing the displacement versus load curves presented in Simo et al \cite{simo1986three}, the beam is discretised into $8$ CVs and an end force of $3000$ N is applied to its free end at $30$ N load increments. Fig. \ref{fig:45degArcLoad} presents the displacement components of the free end of the beam and the results are found to be in close agreement with those reported by Simo et al. \cite{simo1986three}.

To investigate the order of accuracy of the solver in calculating forces, the present test case is reproduced by applying displacements at the right end instead of forces. To obtain mesh independent displacements, a fine mesh of $1280$ CVs is selected and for an applied force of $600$ N along the $z$-axis, the values observed are $w_x = -23.5607$ m, $w_y = -13.6048$ m, $w_z = 53.4756$ m. On applying these displacements to the same test case, the exact force vector $(0 \ 0 \ 600)$ N should be retrieved. Consequently, the test case is run by applying the displacements $(-23.5607, -13.6048, 53.4756)$ m for different mesh sizes, viz. $10$, $20$, $40$, $80$ and $160$ CVs. Fig. \ref{fig:45DegArcMeshConvForce} shows the \% mesh error convergence of the vertical force ($n_z$) for successive mesh size reductions; a quadratic order of error convergence is observed. Since the true value of forces $n_x$ and $n_y$ expected for the given displacements is zero, relative error is not defined; however, the force components obtained for successive mesh reductions are seen to go to zero (Table \ref{table:45DegArcResultsForce}). 

\begin{table}[!]
	\caption{Curved cantilever in $xy-$plane subjected to a force along $z-$axis: force values observed for different mesh sizes} 
	\centering 
	\begin{tabular}{c c c} 
		\hline
		CVs & $n_{x}$ (N) & $n_{y}$ (N) \\ 
		\hline
		$10$ &  $53.98$ & $24.76$ \\
		$20$ &  $12.81$ & $5.91$ \\
		$40$ &  $3.17$ & $1.47$ \\
		$80$ &  $0.80$ & $0.37$ \\
		$160$ &  $0.21$ & $0.09$ \\[0.5ex]
		\hline 
	\end{tabular}
	\label{table:45DegArcResultsForce} 
\end{table}

\subsection{Deep-circular arch with a concentrated in-plane force at the crown location}

\noindent The instabilities of asymmetric clamped-hinged arches were first investigated by DaDeppo and Schmidt \cite{dadeppo1975instability} and then by many other authors \cite{simo1986three,meier2014objective,noor1981mixed}. In the current work, a $215^{\circ}$ deep circular arch with a unit circular cross-section, hinged at the left end ($w_x = w_y = w_z = 0$) and clamped at the right ($w_x = w_y = w_z = 0$ and $\psi_x = \psi_y = \psi_z = 0$), having a radius $R = 100$ m and length $375.245$ m is considered and a force-displacement behaviour of the arch is studied until buckling. A concentrated point load $P \equiv (0 \quad n_y \quad 0)$ acts at the crown location of the arch. The adopted mechanical properties are: $EI_{2} = EI_{3} = GJ = 1 \times 10^4$ Nm$^2$. 

\begin{figure}[!]
	\centering
	\includegraphics[scale = 0.6]{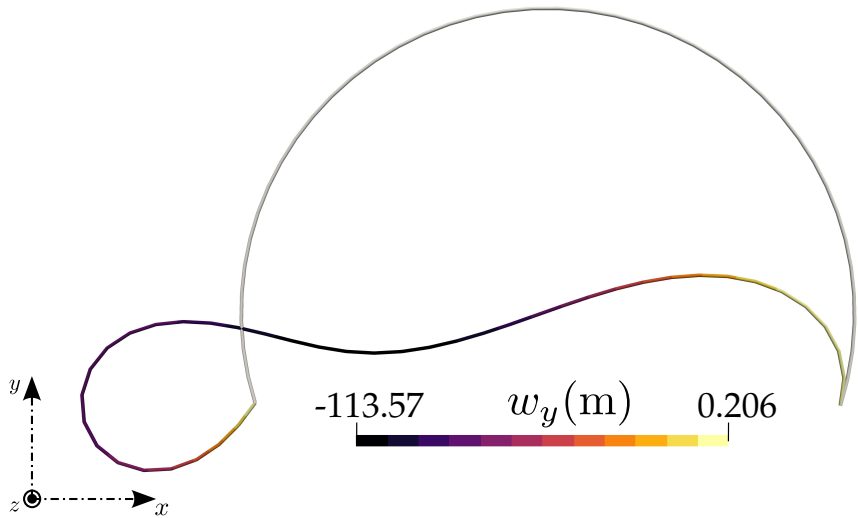}
	\caption{Deep-circular arch with a concentrated in-plane force at the crown location: initial and deformed configuration}
	\label{fig:circularArch}
\end{figure}

\begin{figure}[!]
	\centering
	\includegraphics[width=0.8\textwidth]{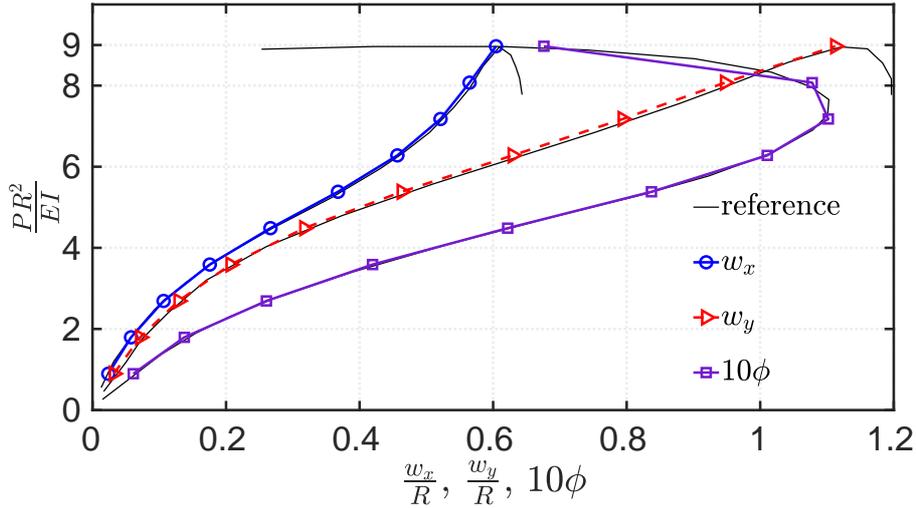}
	\caption{Deep-circular arch with a concentrated in-plane force at the crown location: the normalised horizontal and vertical displacements of the crown point ($w_{x}/R, w_{y}/R$) and the rotation angle ($ 10\phi$) of the crown point for different loading values of the force}
	\label{fig:circularArchResults}
\end{figure}

In accordance with the literature \cite{simo1986three}, the entire spatial domain is discretised into $40$ CVs and Fig. \ref{fig:circularArch} shows the initial and the deformed configuration of the arch. For an increasing load at $1$ N increments up to $8$ N and $0.005$ N increments beyond $8$ N , the predictive critical buckling load is found to be is $9.085$ N exact value where the value reported by DaDeppo and Schmidt \cite{dadeppo1975instability} is $8.97$ N. The solver, being quasi-static, does not converge for the unstable post-buckling analysis of the circular arch; beyond the critical buckling load, a dynamic analysis or special procedures would be required. For the loading values less than $8$ N, the residuals converge at an average of $8$ iterations per pseudo time increment; for loads higher than $8$ N, the residuals converge in an average of $17$ iterations. The total time of computation is less than $1$ s. Fig. \ref{fig:circularArchResults} shows the variation of the dimensionless quantity $\frac{PR^2}{EI}$ versus the normalised displacements ($w_{x}$ and $w_{y}$ with respect to $R$) and rotation angle $\phi$ of the tangent vector at the crown location. The results are found to be in good agreement with the exact solutions presented by DaDeppo and Schmidt \cite{dadeppo1975instability}. 

\section{CONCLUDING REMARKS}
\noindent This paper is the first to develop and verify a total Lagrangian cell-centred finite volume methodology for geometrically exact beams with arbitrary initial curvatures subjected to finite displacements and rotations. The mathematical model and the corresponding FV-based spatial and equation discretisation are described in detail. 
\noindent The potential of the developed methodology has been tested using five complementary benchmark test cases, where the spatial discretisation has been shown to be second-order accurate for displacements and forces. For the cases examined, even though the computational time is quite fast (less than $1$ minute), it is found that the results require a higher number of outer iterations to converge as compared to those reported in literature (average of 6 iterations instead of 2 for pure bending of cantilever beam). Nevertheless, the predictions using the proposed FV methodology are found to be in good agreement with the analytical solutions and the FE based numerical results. Extension of the quasi-static FV formulation to dynamic analysis will be considered in the future work.


\section*{Acknowledgements}
\noindent This publication has emanated from research [conducted with the financial support of/supported in part by a grant from] Science Foundation Ireland under Grant number RC2302\_2. For the purpose of Open Access, the author has applied a CC BY public copyright license to any Author Accepted Manuscript version arising from this submission. Financial support is gratefully acknowledged from the Irish Research Council through the Laureate program, grant number IRCLA/2017/45. Additionally, the authors want to acknowledge project affiliates, Bekaert, through the Bekaert University Technology Centre (UTC) at University College Dublin (\url{www.ucd.ie/bekaert}), and I-Form, funded by Science Foundation Ireland (SFI) Grant Number 16/RC/3872, co-funded under European Regional Development Fund and
by I-Form industry partners. Provision of computational facilities and support from the DJEI/DES/SFI/HEA Irish Centre for High-End Computing (ICHEC, \url{www.ichec.ie}) and ResearchIT Sonic cluster, funded by UCD IT Services and the Research Office is gratefully acknowledged.

%
\section*{Conflict of interest}
\noindent The authors declare that they have no conflict of interest.


\appendix

\setcounter{equation}{0}
\renewcommand{\theequation}{A.\arabic{equation}}

\section{Linearisation of stress-resultants}
\label{app:linearisation}

\noindent The linearisation of the stress-resultants involves a systematic use of the directional derivatives of the kinematic quantities and the strain measures. For a given beam configuration, the linearisation of perturbed kinematic quantities $\big(\bm{r}_\epsilon= \bm{r}_0 + \epsilon \Delta \bm{w},(\RMT)_\epsilon = \exp(\epsilon \widehat{\Delta \RV}) \RMT \big)$ yields,

\begin{equation}
			\frac{d}{d\epsilon} \Big(\bm{r}_\epsilon\Big) \bigg|_{\epsilon = 0} = \Delta \bm{w}
			\qquad \qquad	
		\frac{d}{d\epsilon} \Big((\RMT)_\epsilon\Big) \bigg|_{\epsilon = 0} =  \widehat{ \Delta \RV} \RMT^{\rm *}		
\end{equation}

\noindent where $\epsilon$ is a small scalar perturbation. The linearised forms for the perturbed (material) strain measures $\bm{\Gamma}_\epsilon  = (\RMT^{\rm T})_\epsilon \bm{r}^\prime_\epsilon - \bm{e}_1$ and $\bm{K}_\epsilon = (\RMT^{\rm T})_\epsilon (\RMT^\prime)_\epsilon$ \cite{simo1986three,Marino2016383} are,

\begin{equation}\label{eq:linGamma}
	\frac{d}{d\epsilon} \Big(\bm{\Gamma}_\epsilon \Big) \bigg|_{\epsilon = 0} =(\RMT^{\rm *})^{\rm T} \widehat{(\bm{r}^{*})^\prime} (\Delta \RV) +(\RMT^{\rm *})^{\rm T} \Delta \bm{w}^\prime	
\end{equation}

\begin{equation}
\frac{d}{d\epsilon} \Big( \bm{K}_\epsilon \Big) \bigg|_{\epsilon = 0} =(\RMT^{\rm *})^{\rm T} \Delta \RV^\prime
\end{equation}

\noindent Using the above values, the spatial force, $\bm{n} = \RMT \bm{C}_{\rm N} \bm{\Gamma}$ is linearised as,

\begin{align}
	\begin{split}
		\frac{d}{d\epsilon} \Big( \bm{n}_\epsilon \Big) \bigg|_{\epsilon = 0} &= \frac{d}{d\epsilon} \Big( (\RMT)_\epsilon \bm{C}_{\rm N}  \bm{\Gamma}_\epsilon \Big) \bigg|_{\epsilon = 0} \\[.2cm]
		& = \RMT^{\rm *} \bm{C}_{\rm N} \bm{\Gamma}^{\rm *} +  \frac{d}{d\epsilon} (\RMT)_{\epsilon} \bigg|_{\epsilon = 0} \bm{C}_{\rm N} \bm{\Gamma}^{*} + (\RMT^{\rm *}) \bm{C}_{\rm N} \frac{d}{d \epsilon} (\bm{\Gamma}_\epsilon) \bigg|_{\epsilon=0} \\[.2cm]
		\implies  \bm{n} & = \underbrace{\RMT^{\rm *} \bm{C}_{\rm N} \bm{\Gamma}^{\rm *}}_{\bm{n}^{*}} - \big(\widehat{\RMT^{\rm *} \bm{C}_{\rm N} \bm{\Gamma}^{\rm *}}\big) \Delta \RV + \big(\RMT^{\rm *} \bm{C}_{\rm N}(\RMT^{\rm *})^{\rm T}\big) \Big[ \Delta \bm{w}^\prime + \widehat{(\bm{r}^{*})^\prime} (\Delta \RV) \Big]
	\end{split}
\end{align}

\noindent Similarly, the spatial moment, $\bm{m} = \RMT \bm{C}_{\rm M} \bm{K}$ can be linearised as,

\begin{align}
	\begin{split}
		\frac{d}{d\epsilon} \Big( \bm{m}_\epsilon \Big) \bigg|_{\epsilon = 0} &= \frac{d}{d\epsilon} \Big( (\RMT)_\epsilon \bm{C}_{\rm M}  \bm{K}_\epsilon \Big) \bigg|_{\epsilon = 0} \\[.2cm]
		& = \RMT^{\rm *} \bm{C}_{\rm M} \bm{K}^{\rm *} +   \frac{d}{d\epsilon} (\RMT)_{\epsilon} \bigg|_{\epsilon = 0} \bm{C}_{\rm M} \bm{K}^{*} + (\RMT^{\rm *}) \bm{C}_{\rm M} \frac{d}{d \epsilon} (\bm{K}_\epsilon) \bigg|_{\epsilon=0} \\[.2cm]
		\implies  \bm{m} & = \RMT^{\rm *} \bm{C}_{\rm M} \bm{K}^{\rm *} - \big(\widehat{\RMT^{\rm *} \bm{C}_{\rm N} \bm{K}^{\rm *}}\big) \Delta \RV + (\RMT^{*}\bm{C}_{\rm M}(\RMT^{\rm *})^{\rm T})\Delta\RV^\prime
	\end{split}
\end{align}

\noindent The linearised counterpart of the term $(\bm{r}^\prime\times\bm{n})$ in Eq. (\ref{eq:momentBalance}) is given by,

\begin{align}
	\begin{split}
		\frac{d}{d\epsilon} \Big( \bm{r}^\prime\times\bm{n} \Big) \bigg|_{\epsilon = 0} &= \frac{d}{d\epsilon} \Big( \bm{r}^\prime_\epsilon \times (\RMT)_\epsilon \bm{C}_{\rm N} \bm{\Gamma_\epsilon} \Big) \bigg |_{\epsilon = 0} \\[.2cm]
		& =  (\bm{r}^{*})^\prime \times (\RMT^{*}) \bm{C}_{\rm N} \bm{\Gamma^{*}} + \frac{d}{d \epsilon} (\bm{r}^\prime_\epsilon) \bigg|_{\epsilon = 0} \times \RMT^{*} \bm{C}_{\rm N} \bm{\Gamma}^{*} \\[.2cm]
		& \qquad +(\bm{r}^{*})^\prime \times \frac{d}{d \epsilon} \Big( (\RMT)_\epsilon \Big) \bigg|_{\epsilon=0} \bm{C}_{\rm N} \bm{\Gamma}^{*} + (\bm{r}^{*})^\prime \times (\RMT)^{*}  \bm{C}_{\rm N} \frac{d}{d \epsilon} \Big( \bm{\Gamma}_\epsilon \Big) \bigg|_{\epsilon=0} \\[.2cm]
		\implies \big(\bm{r}^\prime\times\bm{n}\big) &= \widehat{(\bm{r}^{*})^\prime} (\RMT^{*}) \bm{C}_{\rm N} \bm{\Gamma^{*}} + \Delta \bm{w}^\prime \times \big(\RMT^{*} \bm{C}_{\rm N} \bm{\Gamma}^{*}\big) + (\bm{r}^{*})^\prime \times \big(\widehat{ \Delta \RV}\big) \RMT^{\rm *} \bm{C}_{\rm N} \bm{\Gamma^{*}} \\[.2cm]
		& \qquad + (\bm{r}^{*})^\prime \times (\RMT)^{*}  \bm{C}_{\rm N} \Big[\RMT^{\rm *, T} \widehat{(\bm{r}^{*})^\prime} (\Delta \RV) +(\RMT^{\rm *})^{\rm T} \Delta \bm{w}^\prime	\Big] \\[.2cm]
		& = \widehat{(\bm{r}^{*})^\prime} \bm{n}^{*} - \big(\widehat{\bm{n}^{*}}\big) \Delta \bm{w}^\prime +  \widehat{(\bm{r}^{*})^\prime} \big( \bm{n}- \bm{n}^{*}\big)
	\end{split}
\end{align}

\noindent Here, $\big(\cdot\big)^{\rm *}$ denotes the known fields from the previously converged Newton-Raphson iteration. The manipulation of skew-symmetric tensors and the corresponding axial vector relation, $\bm{\theta} \times \bm{h} = \hat{\bm{\theta}}\bm{h} = -(\bm{h} \times \bm{\theta})= - \hat{\bm{h}} \bm{\theta}$ is used in all the above expressions.


\setcounter{equation}{0}
\renewcommand{\theequation}{B.\arabic{equation}}

\section{Coefficients of the block linear system}
\label{app:BlockCoeffs}


\subsection{Force Equilibrium}
\noindent The linearised equation for spatial force $\bm{n}$ (Eq. \ref{eq:linForce}) can be re-written as,

\begin{equation}	
	\bm{n} = \underbrace{\RMT^{*}\bm{C}_{\rm N}\bm{\Gamma}^{*}}_{\bm{C}_{\rm exp \rm w}}
	+ \underbrace{(\RMT^{*} \bm{C}_{\rm N} (\RMT^{\rm *})^{\rm T})}_{\bm{C}_{\rm w \rm w}}\Delta\bm{w}^\prime	
	+ \underbrace{\left[\RMT^{*} \bm{C}_{\rm N} (\RMT^{\rm *})^{\rm T}\widehat{(\bm{r}^{*})^\prime} - \widehat{(\RMT^{*}\bm{C}_{\rm N}\bm{\Gamma}^{*})}\right]}_{\bm{C}_{\rm w \RV}} \Delta \RV
\end{equation}

\noindent Using the discretisation schemes for cell faces (Eq. \ref{eq:faceValue}) and cell-face derivatives (Eq. \ref{eq:faceDerivative}), the discretised force equilibrium (Eq. \ref{eq:discForce}) can be expanded and written in matrix-form as,

\begin{equation}
	\big[\bm{A}_{\rm W} \big]_{\bm{n}}
	\begin{bmatrix}
		(\Delta \bm{w})_{\rm W}\\
		(\Delta \RV)_{\rm W}
	\end{bmatrix}
	+
	\big[\bm{A}_{\rm C} \big]_{\bm{n}}
	\begin{bmatrix}
		(\Delta \bm{w})_{\rm C}\\
		(\Delta \RV)_{\rm C}
	\end{bmatrix}
	+
	\big[\bm{A}_{\rm E} \big]_{\bm{n}}
	\begin{bmatrix}
		(\Delta \bm{w})_{\rm E}\\
		(\Delta \RV){\rm E}
	\end{bmatrix}
	= (\bm{R^{\bm{w}}})_{\rm C}
\end{equation}

\noindent where the coefficient matrices $\big[\bm{A}_{\rm W} \big]_{\bm{n}}$, $\big[\bm{A}_{\rm C} \big]_{\bm{n}}$ and $\big[\bm{A}_{\rm E} \big]_{\bm{n}}$ are $3 \times 6$ matrices containing the contributions of $\rm W$, $\rm C$ and $\rm E$ cell-centres respectively are given by,

\begin{equation}	
	\big[\bm{A}_{\rm W} \big]_{\bm{n}} 
	= 
	\begin{bmatrix} 
		\frac{1}{L_{\rm C}} \Big( \gamma_{\rm w} (\bm{C}_{\rm w \rm w})_{\rm W} + (1 - \gamma_{\rm w})(\bm{C}_{\rm w \rm w})_{\rm C} \Big)
		& \qquad
		-\gamma_{\rm w} \Big( \gamma_{\rm w} (\bm{C}_{\rm w \RV})_{\rm W} + (1 - \gamma_{\rm w}) (\bm{C}_{\rm w \RV})_{\rm C}\Big) \end{bmatrix}	
\end{equation}

\vspace{0.5cm}

\begin{equation}	
	\big[\bm{A}_{\rm C} \big]_{\bm{n}} 
	= 
	\begin{bmatrix} 
	-\frac{1}{L_{\rm C}} \Big( \gamma_{\rm e} (\bm{C}_{\rm w \rm w})_{\rm E} + (1 - \gamma_{\rm e})(\bm{C}_{\rm w \rm w})_{\rm C}
	+ 
	\gamma_{\rm w} (\bm{C}_{\rm w \rm w})_{\rm W} + (1 - \gamma_{\rm w})(\bm{C}_{\rm w \rm w})_{\rm C} \Big)
	\\[0.5cm]
	(1 -\gamma_{\rm e}) \Big( \gamma_{\rm e} (\bm{C}_{\rm w \RV})_{\rm E} + (1 - \gamma_{\rm e}) (\bm{C}_{\rm w \RV})_{\rm C}\Big)
	-
	(1 -\gamma_{\rm w}) \Big( \gamma_{\rm w} (\bm{C}_{\rm w \RV})_{\rm W} + (1 - \gamma_{\rm w}) (\bm{C}_{\rm w \RV})_{\rm C}\Big)  \end{bmatrix}^{T}	
\end{equation}

\vspace{0.5cm}

\begin{equation}	
	\big[\bm{A}_{\rm E} \big]_{\bm{n}} 
	= 
	\begin{bmatrix} 
		\frac{1}{L_{\rm C}} \Big( \gamma_{\rm e} (\bm{C}_{\rm w \rm w})_{\rm E} + (1 - \gamma_{\rm e})(\bm{C}_{\rm w \rm w})_{\rm C} \Big)
		& \qquad
		\gamma_{\rm e} \Big( \gamma_{\rm e} (\bm{C}_{\rm w \RV})_{\rm E} + (1 - \gamma_{\rm e}) (\bm{C}_{\rm w \RV})_{\rm C}\Big) \end{bmatrix}	
\end{equation}

\noindent The residual forces evaluated at the cell centre $\rm C$ is given by,

\begin{equation}
	(\bm{R}^{\bm{w}})_{\rm C}
	= \gamma_{\rm w}(\bm{C}_{\rm exp \rm w})_{\rm W} + (1 -\gamma_{\rm w})(\bm{C}_{\rm exp \rm w})_{\rm C} - \gamma_{\rm e}(\bm{C}_{\rm exp \rm w})_{\rm E} - (1 -\gamma_{\rm e})(\bm{C}_{\rm exp \rm w})_{\rm C} - \bm{f}_{\rm C} L_{\rm C}
\end{equation}

\subsection{Moment Equilibrium}

\noindent The discretised moment equilibrium Eq. \ref{eq:discMoment} can be written as,

\begin{equation}
	\bm{m} = \underbrace{\RMT^{*}\bm{C}_{\rm M}\bm{K}^{*}}_{\bm{C}_{\rm exp \rm m}}
	+ \underbrace{-\widehat{(\RMT^{*} \bm{C}_{\rm M} \bm{K}^{*})}}_{\bm{C}_{\rm m \rm \RV}}\Delta \RV 
	+ \underbrace{(\RMT^{*} \bm{C}_{\rm M}(\RMT^{\rm *})^{\rm T})}_{\bm{C}_{\rm m \rm \RV \rm 2}}\Delta\RV^\prime
\end{equation}

\noindent Using Eq. \ref{eq:linForce}, the terms in the linearised counterpart of the term $(\bm{r}^\prime\times\bm{n})$ in Eq. (\ref{eq:linRcrossN}) can be rearranged and expressed as,

\begin{equation}
	\begin{split}
		(\bm{r}^\prime\times\bm{n}) = \underbrace{\widehat{(\bm{r}^{*})^\prime}(\RMT^{*}\bm{C}_{\rm N}\bm{\Gamma}^{*})}_{\bm{C}_{\rm exp \rm m \rm w}}
		+ \underbrace{\Big[ \widehat{(\bm{r}^{*})^\prime}(\RMT^{*} \bm{C}_{\rm N} (\RMT^{\rm *})^{\rm T})-(\widehat{\RMT^{*}\bm{C}_{\rm N}\bm{\Gamma}^{*}}) \Big]}_{\bm{C}_{\rm m \rm w}}\Delta\bm{w}^\prime
		\\
		+ \underbrace{\Big[ \widehat{(\bm{r}^*)^\prime}(\RMT^{*} \bm{C}_{\rm N} (\RMT^{\rm *})^{\rm T}) \widehat{(\bm{r}^{*})^\prime} -  \widehat{(\bm{r}^*)^\prime} (\widehat{\RMT^{*}\bm{C}_{\rm N}\bm{\Gamma}^{*}}) \Big]}_{\bm{C}_{\rm m \rm w \rm \RV}} \Delta \RV
	\end{split}
\end{equation}

\noindent The discretised moment equilibrium Eq. \ref{eq:discMoment} thus can be expanded and written as,

\begin{equation}
	\big[\bm{A}_{\rm W} \big]_{\bm{m}}
	\begin{bmatrix}
		(\Delta \bm{w})_{\rm W}\\
		(\Delta \RV)_{\rm W}
	\end{bmatrix}
	+
	\big[\bm{A}_{\rm C} \big]_{\bm{m}}
	\begin{bmatrix}
		(\Delta \bm{w})_{\rm C}\\
		(\Delta \RV)_{\rm C}
	\end{bmatrix}
	+
	\big[\bm{A}_{\rm E} \big]_{\bm{m}}
	\begin{bmatrix}
		(\Delta \bm{w})_{\rm E}\\
		(\Delta \RV)_{\rm E}
	\end{bmatrix}
	=
	\begin{bmatrix}
		(\bm{R}^{\bm \RV})_{\rm C}
	\end{bmatrix}
\end{equation}

\noindent The coefficient matrices $\big[\bm{A}_{\rm W} \big]_{\bm{m}}$, $\big[\bm{A}_{\rm C} \big]_{\bm{m}}$ and $\big[\bm{A}_{\rm E} \big]_{\bm{m}}$ are $3 \times 6$ matrices given as,

\begin{equation}	
	\big[\bm{A}_{\rm W} \big]_{\bm{m}}
	= 
	\begin{bmatrix} 
		-\frac{1}{2}\frac{L_{\rm w}}{L_{\rm C}} \Big( \gamma_{\rm w} (\bm{C}_{\rm m \rm w})_{\rm W} + (1 - \gamma_{\rm w})(\bm{C}_{\rm m \rm w})_{\rm C} \Big) \\[0.8cm]
		\begin{gathered}
			-\gamma_{\rm w} \Big( \gamma_{\rm w} (\bm{C}_{\rm m \RV})_{\rm W} + (1 - \gamma_{\rm w}) (\bm{C}_{\rm m \RV})_{\rm C}\Big)
			+ \frac{1}{L_{\rm C}} \Big( \gamma_{\rm w} (\bm{C}_{\rm m \RV 2})_{\rm W} + (1 - \gamma_{\rm w})(\bm{C}_{\rm m \RV 2})_{\rm C} \Big) \\
			+ \frac{1}{2} L_{\rm w} \gamma_{\rm w} \Big( \gamma_{\rm w} (\bm{C}_{\rm m \rm w \RV})_{\rm W} + (1 - \gamma_{\rm w})(\bm{C}_{\rm m \rm w \RV})_{\rm C} \Big) 
		\end{gathered}	
	\end{bmatrix}^{T}
\end{equation}

\begin{equation}	
	\big[\bm{A}_{\rm C} \big]_{\bm{m}} 
	= 
	\begin{bmatrix} 
		\frac{1}{2}\frac{L_{\rm w}}{L_{\rm C}} \Big( \gamma_{\rm w} (\bm{C}_{\rm m \rm w})_{\rm W} + (1 - \gamma_{\rm w})(\bm{C}_{\rm m \rm w})_{\rm C} \Big)
		- \frac{1}{2}\frac{L_{\rm e}}{L_{\rm C}} \Big( \gamma_{\rm e} (\bm{C}_{\rm m \rm w})_{\rm E} - (1 - \gamma_{\rm e})(\bm{C}_{\rm m \rm w})_{\rm C} \Big) \\[0.8cm]
		\begin{gathered}
			(1 -\gamma_{\rm e}) \Big( \gamma_{\rm e} (\bm{C}_{\rm m \RV})_{\rm E} + (1 - \gamma_{\rm e}) (\bm{C}_{\rm m \RV})_{\rm C}\Big)
			- \frac{1}{L_{\rm C}} \Big( \gamma_{\rm e} (\bm{C}_{\rm m \RV 2})_{\rm E} + (1 - \gamma_{\rm e}) (\bm{C}_{\rm m \RV 2})_{\rm C}\Big) \\
			- (1 -\gamma_{\rm w}) \Big( \gamma_{\rm w} (\bm{C}_{\rm m \RV})_{\rm W} + (1 - \gamma_{\rm w}) (\bm{C}_{\rm m \RV})_{\rm C}\Big)
			- \frac{1}{L_{\rm C}} \Big( \gamma_{\rm w} (\bm{C}_{\rm m \RV 2})_{\rm W} + (1 - \gamma_{\rm w}) (\bm{C}_{\rm m \RV 2})_{\rm C}\Big) \\ 
			+ \frac{1}{2} L_{\rm e} (1 -\gamma_{\rm e}) \Big( \gamma_{\rm e} (\bm{C}_{\rm m \rm w \RV})_{\rm E} + (1 - \gamma_{\rm e}) (\bm{C}_{\rm m \rm w \RV})_{\rm C}\Big) \\
			+ \frac{1}{2} L_{\rm w} (1 -\gamma_{\rm w}) \Big( \gamma_{\rm w} (\bm{C}_{\rm m \rm w \RV})_{\rm W} + (1 - \gamma_{\rm w}) (\bm{C}_{\rm m \rm w \RV})_{\rm C}\Big)
		\end{gathered}
	\end{bmatrix}^{T}	
\end{equation}

\vspace{0.5cm}

\begin{equation}	
	\big[\bm{A}_{\rm E} \big]_{\bm{m}} 
	= 
	\begin{bmatrix} 
		\frac{1}{2}\frac{L_{\rm e}}{L_{\rm C}} \Big( \gamma_{\rm e} (\bm{C}_{\rm m \rm w})_{\rm E} + (1 - \gamma_{\rm e})(\bm{C}_{\rm m \rm w})_{\rm C} \Big) \\[0.8cm]
		\begin{gathered}
			\gamma_{\rm e} \Big( \gamma_{\rm e} (\bm{C}_{\rm m \RV})_{\rm E} + (1 - \gamma_{\rm e}) (\bm{C}_{\rm m \RV})_{\rm C}\Big)
			+ \frac{1}{L_{\rm C}} \Big( \gamma_{\rm e} (\bm{C}_{\rm m \RV 2})_{\rm E} + (1 - \gamma_{\rm e})(\bm{C}_{\rm m \RV 2})_{\rm C} \Big) \\
			+ \frac{1}{2} L_{\rm e} \gamma_{\rm e} \Big( \gamma_{\rm e} (\bm{C}_{\rm m \rm w \RV})_{\rm E} + (1 - \gamma_{\rm e})(\bm{C}_{\rm m \rm w \RV})_{\rm C} \Big) 
		\end{gathered}	
	\end{bmatrix}^{T}
\end{equation}

\noindent The residual moment evaluated at the cell centre $\rm C$ is given by,

\begin{equation}
	\begin{gathered}
		(\bm{R}^{\RV})_{\rm C} \quad = \quad \gamma_{\rm w}(\bm{C}_{\rm exp \rm m})_{\rm W} + (1 -\gamma_{\rm w})(\bm{C}_{\rm exp \rm m})_{\rm C} - \gamma_{\rm e}(\bm{C}_{\rm exp \rm m})_{\rm E} - (1 -\gamma_{\rm e})(\bm{C}_{\rm exp \rm m})_{\rm C} \\
		- \frac{1}{2}L_{\rm w} \Big(\gamma_{\rm w}(\bm{C}_{\rm exp \rm m \rm w})_{\rm W} + (1 -\gamma_{\rm w})(\bm{C}_{\rm exp \rm m \rm w})_{\rm C} \Big) \\
		- \frac{1}{2}L_{\rm e} \Big(\gamma_{\rm e}(\bm{C}_{\rm exp \rm m \rm w})_{\rm E} + (1 -\gamma_{\rm e})(\bm{C}_{\rm exp \rm m \rm w})_{\rm C} \Big) - \bm{t}_{\rm C} L_{\rm C}
	\end{gathered} 
\end{equation}


\noindent The diagonal and off-diagonal coefficient matrix contributions from the force and moment equilibrium equations can be assembled and the coupled block linear system (Eq. \ref{eq:discreteEquation}) coefficient matrices take the form as follows,
\begin{equation}	
	\bm{A}_{\rm W}
	= 
	\begin{bmatrix}
		\big[\bm{A}_{\rm W} \big]_{\bm{n}}  \\ \\
		\big[\bm{A}_{\rm W} \big]_{\bm{m}}		
	\end{bmatrix} 	
	\quad ; \quad
	\bm{A}_{\rm C} 
	= 
	\begin{bmatrix}
		\big[\bm{A}_{\rm C} \big]_{\bm{n}}  \\ \\
		\big[\bm{A}_{\rm C} \big]_{\bm{m}}			
	\end{bmatrix}
	\quad ; \quad
	\bm{A}_{\rm E} 
	= 
	\begin{bmatrix}
		\big[\bm{A}_{\rm E} \big]_{\bm{n}}  \\ \\
		\big[\bm{A}_{\rm E} \big]_{\bm{m}}			
	\end{bmatrix}
\end{equation}

\bibliographystyle{elsarticle-harv}

\bibliography{references}   

\end{document}